\documentclass[12pt]{article}
\usepackage{amsfonts,amsthm,amssymb,amsmath}
\title{Genome of Descartes Folium via Normalization}
\author{\framebox[1.1\width]{Adrian Constantinescu}, Constantin Udri\c ste, Stelu\c ta Pricopie}
\date{}
\usepackage{amssymb}
\usepackage{eucal}
\usepackage{amsmath}

\begin{document}
\maketitle
\newtheorem{Th}{Theorem}
\newtheorem{Co}{Corollary}
\newtheorem{Prop}{Proposition}
\newtheorem{Lem}{Lemma}
\newtheorem{Def}{Definition}

\begin{abstract}
The Folium of Descartes in $\mathbb{K}\times\mathbb{K}$ carries group laws,
defined entirely in terms of algebraic operations over the field $\mathbb{K}$.
The problems discussed in this paper include:
normalization of Descartes Folium,
group laws and morphisms, exotic structures, exotic structures,
second exotic structure, some topologies on Descartes Folium, differential structure on Descartes Folium,
first isomorphism of algebraic Lie groups over $\mathbb{K}$,
second isomorphism of algebraic Lie groups over $\mathbb{K}$, derived structures
of algebraic Lie groups, a differential/complex analytic structure on Descartes Folium,
Descartes Folium as a topological field, etc.
For predicting these terms, we focus on methods that exploit diagram manipulation techniques
(as alternatives to algebraic method of proofs).
All our results confirm that the Descartes Folium stores
natural group structures, unsuspected till now.
\end{abstract}



\section{Group Structure on Descartes Folium}

\hspace{0.5cm} Traditionally, the group laws were analyzed on regular elliptic curves (\cite{[1]}, \cite{[5]}-\cite{[13]}, \cite{[16]}).
Our theory refers to the {\it Descartes Folium} which is a non-smooth curve
$$DF: x^3 + y^3 -3axy=0 \subset A^2_\mathbb{K},\, a \in \mathbb{K}\setminus\{0\}$$
and to its {\it projective closure} defined by homogenization, i.e.,
$$\overline{DF}: x^3 + y^3 -3axyz=0 \subset P^2_\mathbb{K}$$
and called the {\it projective  Descartes Folium}, too.

The Descartes Folium $\overline{DF}$
is a non-smooth cubic (with a singular point, $O=(0,0)$), non-isomorphic with an elliptic curve,
that admit a multiplicative group structure (see [11]). Now wide research area, highlighting the
group structures by means of canonical isomorphisms. The phrase "Genome of Descartes Folium" means
"Group Structures on Descartes Folium".

The description of a group law on the Descartes Folium can be summed up as follows:

{\bf Theorem} {\it A group law on $\overline{DF}$ is determined
by a choice of an identity point $I \in \overline{DF}\setminus\{(0,0)\}$ and declaring that if
three points $P, Q, R \in \overline{DF}$ lie on the same straight line (counted with
multiplicity), then $P \star  Q \star R = I$.}

We will always choose the "point at infinity" as the
identity $I$ of the Descartes Folium an then it is not
necessary to eliminate the critical point $(0,0)$. Always, $P^3=I$.

A similar Theorem is well-known for the group structures on elliptic curves (\cite{[16]}).

The study of Descartes Folium is a fascinating subject (\cite{[14]}). With the definition of a
curious group law, Descartes Folium becomes powerful computational devices in number
theory. Perhaps more interesting is that, through careful construction of Descartes Folium,
one can create curves whose group law is identical to that of multiplication
or addition. In a sense, all the operations we use in day to day life can be created
and studied on Descartes Folium, as can some far more exotic ones.

The first oral conjecture about the existence of a group structure on Descartes Folium
was made by C. Udriste, in the intention to provide a pertinent topic for a Doctoral Thesis 
of his PhD student, S. Pricopie, who insisted for new algebraic structures on plane curves. 
The affirmative solution was presented in \cite{[14]}. A. Constantinescu was joined our research group 
being conquered by the novelty and complexity of our research subject \cite{[3]}. 

To develop our theory, we use mainly two fields $\mathbb{K} = \mathbb{R}, \mathbb{C}$.
Also we denote
$A_\mathbb{R} = \mathbb{R}[x, y]$ or $A_\mathbb{C} = \mathbb{C}[x, y]$.
Generically, $A_\mathbb{K}=\mathbb{K}[x,y]$.

\section{Normalization of Descartes Folium}

\hspace{0.5cm} Let ${\cal C}$ be an algebraic curve over an algebraically closed field $\mathbb{K}$. According to the general theory of
normalization of an algebraic variety, by the normalization of an algebraic curve ${\cal C}$ we mean a pair
$(\tilde{\cal C}, \tilde p)$, where

1) ${\cal C}$ is a nonsingular algebraic curve over $\mathbb{K}$;

2) $\tilde p: \tilde{\cal C} \to {\cal C}$ is a finite surjective birational morphism.

(Recall that $\tilde p$ is finite if and only if it is proper and all its fibers are finite.
The morphism $\tilde p$ is birational if and only if there exist Zariski-open subsets $\tilde U\subseteq \tilde{\cal C}$
and  $U\subseteq {\cal C}$ such that $\tilde p(\tilde U) = U$ and $\tilde p: \tilde U {\stackrel{\sim}{\longrightarrow}}U$
is an isomorphism of algebraic varieties).

The pair $(\tilde{\cal C}, \tilde p)$ is uniquely determined up to an isomorphism, i.e., if the pair
$\left({\tilde{\tilde{\cal C}}}, \tilde{\tilde p}\,\right)$
is another normalization on ${\cal C}$, then there exists an isomorphism
$\varphi: {\tilde{\cal C}}{\stackrel{\sim}{\longrightarrow}}\tilde{\tilde{\cal C}}$
of algebraic $\mathbb{K}$-varieties such that $\tilde p = \tilde{\tilde p}\circ \varphi$.

In the sequel we describe the normalizations of the Descartes Folium
$\overline{DF}$, resp. $DF$, over an algebraically closed field $\mathbb{K}$, via some natural parametrizations.

\subsection{Normalization of projective Descartes Folium \\(Parametrization 1)}

\hspace{0.5cm} Let $\mathbb{K}$ be a field with $\hbox{char.}\,\mathbb{K}\not=3$ and
$$\overline{DF}\subset \mathbb{P}_{\mathbb{K}}^2:\,\, x^3 + y^3 - 3axyz=0,\,\,a \in \mathbb{K}\setminus\{0\}.$$
The points at infinity of $\overline{DF}$ are given by the equations $z=0$ and $x^3+y^3=0$.

The equation $\lambda^3+1=0$ has $1$ or $3$ distinct solutions in $\mathbb{K}$ (depending on the field $\mathbb{K}$):
$$\hbox{either}\,\lambda_1=-1\, \hbox{or} \,\lambda_1=-1, \lambda_2=\epsilon_1, \lambda_3=\epsilon_2,\,
\hbox{with}\,\epsilon_i^2-\epsilon_i+1=0, i=1,2.$$ These correspond to $1$ or $3$ distinct points
of $\overline{DF}$ at infinity:
$$\hbox{either}\,\,I=(1,-1,0)\,\,\hbox{or}\,\,I=(1,-1,0), \infty_1=(1,\epsilon_1,0), \infty_2=(1,\epsilon_2,0).$$

Let us consider the following parametrization of $\overline{DF}$:
$$
\begin{array}{ccc} \overline{DF} & (x=3at, y = 3at^2, z=1+t^3)&(x,y,z)\in \overline{DF}\setminus\{O\}\\ \
\bar{p}\Big\uparrow& \Big\uparrow&\Big\downarrow\\ \
P^1_\mathbb{K}=A^1_\mathbb{K}\cup\{\infty\}&t\in A^1_\mathbb{K}&t=\frac{y}{x}\,\,,\end{array}
$$
where we indicated the definition of $\bar p$ and of a partial inverse of $\bar p$.
We complete this diagram by the correspondences
$$
\begin{array}{ccccc}O=(0,0,1)& O=(0,0,1) & I=(1,-1,0)&\infty_1=(1,\epsilon_1,0)&\infty_2=(1,\epsilon_2,0)\\ \
\Big\uparrow&\Big\uparrow& \Big\uparrow&\Big\uparrow&\Big\uparrow\\ \
t=\infty&t=0&t=-1&t=\epsilon_1&t=\epsilon_2\,\,,\end{array}
$$
which spotlights all possible three points at infinity. Hence $\bar p(0)=\bar p(\infty)=O$
and $\bar p|_{\mathbb{P}^1_{\mathbb{K}}}\setminus\{0,\infty\}$ is injective.

If $\mathbb{K}$ is algebraically closed, then $\bar p$ is a proper morphism
since $\mathbb{P}^1_{\mathbb{K}}$ and $\overline{DF}$ are projective algebraic varieties (\cite{[7]}).
 It is easy to verify that $\bar p$ has finite fibers, is surjective
and birational. Consequently the pair $(\mathbb{P}^1_{\mathbb{K}}, \bar p)$ is a normalization of $\overline{DF}$.

Let us mention that the definition of $\bar p$ is of pure natural geometric nature:
for each $t\in \mathbb{P}^1_{\mathbb{K}}, t\not =0, \infty,$ $p(t)$ is the intersection
point in $\mathbb{P}^2_{\mathbb{K}}$ of $\overline{DF}$, different of $O=(0,0,1)$,
with the Zariski projective closure of the affine straight line $y=tx$.

\subsection{Normalization of affine Descartes Folium \\(Parametrization 2)}

\hspace{0.5cm} Let $\mathbb{K}$ be a field with $\hbox{char}.\,\mathbb{K}\not = 3$ and
$$DF\subset \mathbb{A}_{\mathbb{K}}^2:\,\, x^3 + y^3 - 3axyz=0,\,\,a \in \mathbb{K}\setminus\{0\}.$$

We have $DF=\overline{DF}\setminus\{I\}$ or $DF=\overline{DF}\setminus\{I,\infty_1,\infty_2\}$
(depending on the field $\mathbb{K}$) and
$${\bar p}\,^{-1}(\{I\})=\{-1\},\,{\bar p}\,^{-1}(\{\infty_i\})=\{\epsilon_i\},\,i=1,2.$$

Let $M = \mathbb{P}^1_\mathbb{K}\setminus\{t\in \mathbb{A}^1_{\mathbb{K}}|t^3+1=0\}$.
From the parametrization $\bar p$ above, it follows the following parametrization of $DF$:
$$
\begin{array}{cccc} DF & \left(x=\frac{3at}{1+t^3}, y = \frac{3at^2}{1+t^3}\right)&O=(0,0)&(x,y)\in DF\setminus\{O\}\\ \
p\Big \uparrow & \Big \uparrow&\Big \uparrow &\Big \downarrow \\ \
M&t\in \mathbb{A}^1_{\mathbb{K}}&t=\infty&t=\frac{y}{x}\,.\end{array}
$$
We have $p(0)=O$ as well $p(\infty)=O$.

In the situation when $\mathbb{K}$ is algebraically closed,
since $\mathbb{P}^1_\mathbb{K}\setminus\{t\in \mathbb{A}^1_{\mathbb{K}}|t^3+1=0\}={\bar p}\,^{-1}(DF)$
and $\bar p$ is a proper morphism, it follows that $p$ is also proper, has finite fibers, is surjective and birational.
Consequently the pair
$(\mathbb{P}^1_\mathbb{K}\setminus\{t\in \mathbb{A}^1_{\mathbb{K}}|t^3+1=0\},p)$
is a normalization of $DF$.

As for previous parametrization $1$, for $\{t\in \mathbb{A}^1_{\mathbb{K}}|t^3+1=0\}$, $t\not =0$, $p(t)$
is the intersection in $\mathbb{A}^2_{\mathbb{K}}$, different of $O$, of $DF$ with the straight line $y=tx$.
In particular, if $\mathbb{K}=\mathbb{R}$ (or, more general, if $\mathbb{K}$ is such that the equation
$\lambda^3+1=0$ has only the solution $\lambda =-1$ in $\mathbb{K}$), the previous morphism becomes
$$
\begin{array}{cccc} DF & \left(x=\frac{3at}{1+t^3}, y = \frac{3at^2}{1+t^3}\right)&O=(0,0)&(x,y)\in DF\setminus\{O\}\\ \
p\Big \uparrow & \Big \uparrow&\Big \uparrow &\Big \downarrow \\ \
N&t\in \mathbb{A}^1_{\mathbb{K}}&t=\infty&t=\frac{y}{x}\,,\end{array}
$$
where $N = (\mathbb{A}^1_{\mathbb{K}}\setminus\{-1\})\cup\{\infty\}$.

\section{First "Exotic" Structure}

\hspace{0.5cm} Let
$\mathbb{K}=\mathbb{R}$, $DF \subset \mathbb{A}^2_\mathbb{R}$, $\mathbb{A}^1_\mathbb{R}\setminus\{-1\}=\mathbb{R}\setminus\{-1\}$.
The diagram
$$
\begin{array}{cccc} DF & \left(x=\frac{3at}{1+t^3}, y = \frac{3at^2}{1+t^3}\right)&(x,y)\in DF\setminus\{(0,0)\}&(0,0)\\ \
{\wr}\Big\uparrow p& \Big\uparrow&\Big\downarrow&\Big\downarrow\\ \
\mathbb{R}\setminus\{-1\}&t&t=\frac{y}{x}&t=0\end{array}
$$
shows that the function $p$ is bijective. On the other hand
the function $t\to \tau = t+1$, with the inverse $t=\tau - 1{\stackrel{\alpha}{\longleftarrow}} \tau$, is a bijection
$$\mathbb{R}\setminus\{-1\}\underset{\sim}{\stackrel{\alpha}{\longleftarrow}} \mathbb{R}\setminus \{0\}.$$

It appears the diagram
$$
\begin{array}{cccc} DF & \left(x=\frac{3a(\tau -1)}{1+(\tau - 1)^3}, y = \frac{3a(\tau-1)^2}{1+(\tau -1)^3}\right)&(x,y)\in DF\setminus\{(0,0)\}&(0,0)\\ \
{\wr}\Big\uparrow p\,\alpha& \Big\uparrow&\Big\downarrow&\Big\downarrow\\ \
R\setminus\{0\}&\tau &\tau =\frac{y}{x}+1&1\,\end{array}
$$
which proves that $p\,\alpha$ is a bijection. It follows that the group structure on $R\setminus\{0\}$
transfers to $DF$:
$$(p\alpha)(\tau)\circ (p\alpha)(\tau^{\prime})\,{\stackrel{def}{=}}\,(p\alpha)(\tau\tau^{\prime})$$
or
$$
\left(\frac{3a(\tau -1)}{1+(\tau - 1)^3},\frac{3a(\tau -1)^2}{1+(\tau - 1)^3}\right)\circ
\left(\frac{3a(\tau^{\prime} -1)}{1+(\tau^{\prime} - 1)^3},\frac{3a(\tau^{\prime} -1)^2}{1+(\tau^{\prime} - 1)^3}\right)$$
$$\,{\stackrel{def}{=}}\, \left(\frac{3a(\tau\tau^{\prime} -1)}{1+(\tau\tau^{\prime} - 1)^3},\frac{3a(\tau\tau^{\prime} -1)^2}{1+(\tau\tau^{\prime} - 1)^3}\right).$$

It appears an isomorphism of groups
$$(DF,\circ)\simeq (\mathbb{R}\setminus\{0\},\cdot)$$
$$(p\,\alpha)(\tau)\leftarrow \tau.$$
Here, $O=(0,0)\in DF$ is the neutral element for $(DF,\circ)$.

{\bf Remark} More generally, we can define the previous group structure $(DF,\circ)$ and the group
isomorphism $p\alpha$, if we consider a base field $\mathbb{K}$ with $\hbox{char.}\,\mathbb{K} \not = 3$
and with the property that the equation $\lambda^3 + 1=0$ has in $\mathbb{K}$ only the root $\lambda =-1$.
\vspace{0.3cm}

Let us return to the particular case $\mathbb{K}= \mathbb{R}$.
The adjective "exotic" refers to the following explanations. First let us point out that $(DF,\circ)$
is not a topological group, where $DF\subset \mathbb{A}_{\mathbb{R}}^2=\mathbb{R}^2$ is naturally
endowed with the topology induced by the natural real topology of $\mathbb{A}_{\mathbb{R}}^2=\mathbb{R}^2$.
For this let us formulate the following

{\bf Proposition} {\it Let $\mathbb{K}=\mathbb{R}, \mathbb{C}$ and $X=DF\subset \mathbb{A}_{\mathbb{K}}^2$
or $X=\overline{DF}\subset\mathbb{P}_{\mathbb{K}}^2$, with the topology induced by the natural real,
respectively complex, topology of $\mathbb{A}_{\mathbb{K}}^2$, respectively of $\mathbb{P}_{\mathbb{K}}^2$.
Then the topological space $(X,\tau)$ does not admit a topological group structure.}

{\bf Proof} By contrary, let us suppose that there exists a topological group law, denoted by $\cdot$, on $(X,\tau)$.
Let $O\in X$ be the singular point of $X$ (i.e., $O=(0,0)$ or $O=(0,0,1)$) and $O^\prime\in X,\,O^\prime \not=O$.
Then $O^\prime$ is a nonsingular point of $X$.

Let $\theta: X\stackrel{\sim}\longrightarrow X$ be the translation map (which is a homeomorphism) such that
$\theta(O)= O^\prime$ (i.e., $\theta(P)=P\cdot O^{-1}\cdot O^\prime$, for each $P\in X$). Let ${\cal V}$, resp ${\cal V}^\prime$,
be the set of all open neighborhoods of the points $O$, resp. $O^\prime$.

If $\mathbb{K}=\mathbb{R}$, then there exists a decreasing fundamental system of
neighborhoods $\{V_n\}_{n\in \mathbb{N}}\subset{\cal V}$,
resp. $\{V_n^\prime\}_{n\in \mathbb{N}}\subset{\cal V}^\prime$, such that the topological space
$V_n\setminus\{O\}$, resp. $V_n^\prime\setminus\{O^\prime\}$, has $4$, resp. $2$, connected components. Then
$$\underset{V\in {\cal V}}{\underset{\longrightarrow}{\lim}}H^0(V,\mathbb{Z})=\underset{n\in \mathbb{N}}{\lim}H^0(V_n,{\mathbb Z})\approx \mathbb{Z}^4$$
and
$$\underset{V^\prime\in {\cal V}}{\underset{\longrightarrow}{\lim}}H^0(V^\prime,\mathbb{Z})=\underset{n\in \mathbb{N}}{\lim}H^0(V^\prime_n,{\mathbb Z})\approx \mathbb{Z}^2.$$

If $\mathbb{K}=\mathbb{C}$, we have a similar situation, with
$V_n\setminus\{O\}$, resp. $V_n^\prime\setminus\{O^\prime\}$, having $2$, resp. $1$, connected components. Then
$$\underset{V\in {\cal V}}{\underset{\longrightarrow}{\lim}}H^0(V,\mathbb{Z})\approx \mathbb{Z}^2$$
and
$$\underset{V^\prime\in {\cal V}}{\underset{\longrightarrow}{\lim}}H^0(V^\prime,\mathbb{Z})\approx \mathbb{Z}.$$
Consequently, in both cases, $\underset{V\in {\cal V}}{\underset{\longrightarrow}{\lim}}H^0(V,\mathbb{Z})$ is not isomorphic to
$\underset{V^\prime\in {\cal V}}{\underset{\longrightarrow}{\lim}}H^0(V^\prime,\mathbb{Z})$,
which contradicts that $\theta$ is a homeomorphism.

\section{Branches, topologies and \\differential structure}

\hspace{0.5cm} Suppose $\mathbb{K}=\mathbb{R}$. Let us consider the branches of the singularity $O=(0,0)$ of $DF$ as follows:

- the "South branch" $S=p(-1,1)$;

- the "West branch" $W= p(1,\infty)\cup\{O\}\cup p(-\infty,-1)$.

We have $S\cap W=\{O\}$ and $S\cup W= DF \setminus \{V\}$, where $V=p(1)= \left(\frac{3}{2},\frac{3}{2}\right)$ is the "vertex" of $DF$.
The branches $S$ and $W$ are symmetric w.r.t. the first bisector $x-y=0$ of $\mathbb{A}^2_{\mathbb{R}}$.
This means that applying the symmetry $\sigma$ w.r.t. the bisector $x-y=0$, given by
$$
\begin{array}{ccc}\mathbb{A}^2_{\mathbb{R}}&{\stackrel{\sigma}{\stackrel{\sim}{\longrightarrow}}}&\mathbb{A}^2_{\mathbb{R}}\\ \
(x,y)& \longrightarrow &(y,x),\end{array}
$$
we have $\sigma(DF)=DF$ and the branches $S$ and $W$ interchange by $\sigma$
(i.e., $\sigma(S)=W$ and $\sigma(W)=S$).

By the parametrization on $p$, the point $O=p(0)$ is reached on the branch $S$
and it is not reached on the branch $W$.

Let us consider the parametrization on $p^{\prime}=\sigma p$ of $DF$. Then, by using the
interchange of $S$ and $W$ by $\sigma$, it follows that $W=p^{\prime}(-1,1)$ and
$S= p^{\prime}(1,\infty)\cup\{O\}\cup p^{\prime}(-\infty,-1)$.

By the parametrization $p^{\prime}=\sigma p$, the point $\{O\}$ is reached only on the branch
$W$ (not on the branch $S$).

It is easy to see that the pair
$((\mathbb{A}^1_{\mathbb{R}}\setminus\{-1\})\cup\{\infty\}, p^{\prime}=\sigma p)$
is also a normalization of $DF$ (see $\S$ 2).

\subsection{Some topologies on affine Descartes Folium}

\hspace{0.5cm} Concerning the topological properties of the map $p$, we have the following

{\bf Proposition} {\it Suppose $\mathbb{K}=\mathbb{R}$ and $DF\subset \mathbb{A}^2_{\mathbb{R}}$
endowed with the topology $\tau$ induced by the real topology of $\mathbb{A}^2_{\mathbb{R}}$.
Then

(i) the bijective map $p: \mathbb{R}\setminus\{-1\}\stackrel{\sim}\longrightarrow  DF$
is continuous but not a homeomorphism;

(ii) $p|_{\mathbb{R}\setminus\{0,-1\}}:\mathbb{R}\setminus\{0,-1\}\stackrel{\sim}\longrightarrow DF\setminus\{O\}$
is a homeomorphism.

Similar properties hold for the map $p^{\prime}=\sigma p: \mathbb{R}\setminus\{-1\}\stackrel{\sim}\longrightarrow DF$.}

{\bf Proof} (i) Suppose, by contrary, that
$p: \mathbb{R}\setminus\{-1\}\stackrel{\sim}\longrightarrow  DF$
is a homeomorphism. Then
$p\alpha : \mathbb{R}\setminus\{0\}\stackrel{\sim}\longrightarrow  DF$
is also a homeomorphism. Since $(\mathbb{R}\setminus\{0\}, +)$ is a topological group
and $p\alpha$ is a group isomorphism onto $(DF,\circ)$, it follows easy that $(DF,\circ)$
is a topological group, which is not possible.

An alternative proof based on the different connection properties of $\mathbb{R}\setminus\{0\}$
and $DF$ can be done.

(ii) The inverse map
$$
\begin{array}{ccc}\mathbb{R}\setminus\{0,-1\}&{\stackrel{p^{-1}}{\stackrel{\sim}{\longleftarrow}}}&DF\setminus\{O\}\\ \
t=\frac{y}{x}& \longleftarrow &(x,y)\end{array}
$$
is also continuous.

In the previous Proposition we have worked with
the topology $\tau$ on $DF$ which is induced on $DF\subset \mathbb{A}^2_{\mathbb{R}}$
by the real topology of $\mathbb{A}^2_{\mathbb{R}}$. Now let us change the topology $\tau$ on $DF$
with the topology $\tau_S$ (resp. $\tau_W$) defined as follows:

{\bf Definition} {\it $\tau_S$ (resp. $\tau_W$) is the image on $DF$ of the real topology of
$\mathbb{R}\setminus\{-1\}$ by the bijective map
$p: \mathbb{R}\setminus\{-1\}\stackrel{\sim}\longrightarrow  DF$
(resp. by $p^{\prime}=\sigma p: \mathbb{R}\setminus\{-1\}\stackrel{\sim}\longrightarrow  DF$)}

Hence the new topology $\tau_S$ (resp. $\tau_W$) on $DF$ is obtained by carrying the real topology
of an open subset of the normalization of $DF$ by the normalization map $p$ (resp. $p^{\prime}$). It follows that
the topology $\tau_S$ (resp. $\tau_W$) is separated, paracompact and locally compact,  and with countable basis, as well as the fact that
$S= p(-1,1)$ (resp. $W = p^\prime(-1,1)$) is open in $DF$ w.r.t. $\tau_S$ (resp. $\tau_W$). Moreover,
the topological space $(DF,\tau_S)$ (resp. $(DF,\tau_W)$) has two connected components.

\subsection{Some properties of $\tau_S$ (resp. $\tau_W$)}

\hspace{0.5cm} (i) $\tau_S$ (resp. $\tau_W$) is a finer topology than $\tau$ (i.e., $\tau_S, \tau_W\succ \tau$).

(ii) The induced topology $\tau_S|_{DF\setminus\{O\}}$ ($\tau_W|_{DF\setminus\{O\}}$) on
$DF\setminus\{O\}\subset \mathbb{A}^2_{\mathbb{R}}$ is that induced on $DF\setminus\{O\}$ by the real topology
of $\mathbb{A}^2_{\mathbb{R}}$.

Equivalently,
$$\tau_S|_{DF\setminus\{O\}}= \tau_{DF\setminus\{O\}}\,\, \left(\hbox{resp.}\,\tau_W|_{DF\setminus\{O\}}= \tau_{DF\setminus\{O\}}\right).$$

(iii) If $\{U_n\}_{n\in \mathbb{N}}$ is a fundamental system of open neighborhoods of $O$ in
$\mathbb{A}^2_{\mathbb{R}}$, with respect to the real topology, then $\{U_n\cap S\}_{n\in \mathbb{N}}$
(resp. $\{U_n\cap W\}_{n\in \mathbb{N}}$) is a fundamental system of open neighborhoods of $O\in DF$,
in $DF$, with respect to the topology $\tau_S$ (resp. $\tau_W$).

(iv)
$$\tau\cup\{U\cap S\,|\, U\subseteq \mathbb{A}^2_{\mathbb{R}} \,\,\hbox{open subset}\}$$
$$\left(\hbox{resp.}\,\tau\cup\{U\cap W\,|\, U\subseteq \mathbb{A}^2_{\mathbb{R}}\, \,\hbox{open subset}\}\right)$$
is a basis for the topology $\tau_S$ (resp. $\tau_W$). Moreover, for each $V\in \tau_S$ (resp. $\tau_W$),
$$V=(U^\prime \cap DF)\cup(U\cap S)\,\,\, (\hbox{resp.}\,\, V=(U^\prime\cap DF)\cup(U\cap W)),$$
with $U^\prime, U\subseteq \mathbb{A}^2_{\mathbb{R}}$ open subsets.

(v) Let
$$\pi: DF \to \mathbb{R}, \pi(x,y) = \left\{\begin{array}{ccc} \frac{y}{x}&if&(x,y)\not=O\\ \
0&if &(x,y)=O\end{array}\right.$$
$$\left(\hbox{resp.}\,\,\pi^{\prime}: DF \to \mathbb{R}, \pi^{\prime}(x,y) = \left\{\begin{array}{ccc} \frac{x}{y}&if&(x,y)\not=O\\ \
0&if &(x,y)=O\end{array}\right.\right).$$
Then $\tau_S$ (resp. $\tau_W$) is the weakest topology on $DF$ such that $\pi$ (resp. $\pi^{\prime}$)
is continuous ($\mathbb{R}$ endowed with the real topology).

(vi) $\{O\}\subset W$ (resp. $\{O\}\subset S$ is a connected component of the subspace $W$ (resp. $S$)
w.r.t. the topology $\tau_S$ (resp. $\tau_W$). Moreover
$$W=p(1,\infty)\cup\{O\}\cup p(-\infty,-1)$$
$$\left(\hbox{resp.}\,\,S=p^{\prime}(1,\infty)\cup\{O\}\cup p^{\prime}(-\infty,-1)\right)$$
is the representation of $W$ (resp. $S$) as the union of its connected components w.r.t. $\tau_S$ (resp. $\tau_W$).
On the other hand $S$ (resp. $W$) is connected w.r.t. $\tau_S$ (resp. $\tau_W$).

{\bf Proof} Properties (i), (ii) and (v) are direct consequences of the definition of $\tau_S$ (resp. $\tau_W$)
and of the fact that the maps
$$p, p^{\prime}: \mathbb{R}\setminus\{-1\}\stackrel{\sim}\longrightarrow  DF$$
( having $\pi, \pi^\prime: DF \to \mathbb{R}\setminus\{-1\}$ as inverse maps) are continuous and
$$p, p^{\prime}: \mathbb{R}\setminus\{0,-1\}\stackrel{\sim}\longrightarrow  DF\setminus\{O\}$$
are homeomorphisms, where $DF$ (resp. $DF\setminus\{O\}$) above is endowed with the topology $\tau$
(resp. $\tau|_{DF\setminus\{O\}}$).

For property (iii), let us point out firstly that $U\cap S$ (resp. $U\cap W$) is an open subset of $DF$ w.r.t.
$\tau_S$ (resp. $\tau_W$), in particular an open neighborhood of the point $O\in DF$ w.r.t. $\tau_S$
(resp. $\tau_W$), if $U\subseteq \mathbb{A}^2_{\mathbb{R}}$ is an open subset w.r.t.
$\tau_S$ (resp. $\tau_W$), resp. an open neighborhood of the point $O$ in $\mathbb{A}^2_{\mathbb{R}}$.
In fact, $S$ (resp. $W$) is open in $DF$ w.r.t. $\tau_S$ (resp. $\tau_W$) and $U\cap DF \in \tau \subseteq\tau_S$
(resp. $U\cap DF\in \tau_W$) and so $U\cap S = (U\cap DF)\cap S \in \tau_S$ (resp. $\tau_W$).

To end the proof of (iii) it suffices to resume to the topology $\tau_S$ and to show that for an open
neighborhood $V$ of $O$ in $DF$ w.r.t $\tau_S$, there exists an open neighborhood $U$ of $O$ in
$\mathbb{A}^2_{\mathbb{R}}$ such that $V\supseteq U\cap S$. Indeed, we can reduce the situation to
the case $V=p((-\delta,\delta))$, with $0<\delta<0$, because always for such $V$ we have $V\supseteq p((-\delta,\delta))$,
with $0<\delta<1$, and $p((-\delta,\delta))$ is an open neighborhood of $O$ in $DF$ w.r.t. $\tau_S$.

For $t \in p((-\delta,\delta))$, with $0<\delta<1$, we have $|t|<\delta<1$ and from the relation
$x=\frac{3at}{1+t^3}$, where $t=\pi(x,y)$, with $(x,y)\in S$, it follows $3at=(1+t^3)x$ and
$$3|a|\,|t|\leq |1+t^3|\,|x|\leq (1+|t|^3)|x| < (1+\delta^3)|x|< 2|x|.$$
Hence $|t|< \frac{2}{3|a|}\,|x|$. If we consider $\epsilon$, with $0<\epsilon <\delta$, and
$$U = \{(x,y)\in \mathbb{A}^2_{\mathbb{R}}= \mathbb{R}^2\Big|x| <\frac{3|a|}{2}\,\epsilon\},$$
then $U\subseteq  \mathbb{A}^2_{\mathbb{R}}$ is open w.r.t. the standard real topology and
we have $|t|<\epsilon$, for each $(x,y)\in U\cap S$, i.e., $|\pi(x,y)|<\epsilon$, for $(x,y)\in  U\cap S$.
Therefore $\pi(U\cap S)\subset (-\epsilon,\epsilon) \subset (-\delta,\delta)$. Since $\pi = p^{-1}$,
we have then
$$U\cap S = p\pi(U\cap S)\subset p(-\epsilon, \epsilon)\subset p(-\delta,\delta)=V.$$

For property (iv), recall firstly that
$$\tau = \{U\cap DF\,\,|\,\,U \subseteq \mathbb{A}^2_{\mathbb{R}}\,\,\hbox{open subset}\}.$$
Also, we resume to the topology $\tau_S$. Then the family
$$\tau \cup \{U\cap S\,\,|\,\,U \subseteq \mathbb{A}^2_{\mathbb{R}}\,\,\hbox{open subset}\}$$
is closed w.r.t. the finite intersections.

Let $V\subseteq DF$ be an open subset w.r.t. $\tau_S$. If $O\in V$, then $V\supseteq U\cap S$
with $U\subseteq \mathbb{A}^2_{\mathbb{R}}$ open and $O\in U$, according to (iii) and its proof.
If $P\in V$, $P\not = O$, then
$$V\supseteq V\cap (DF\setminus\{O\})\in \tau$$
according to (ii) and $P\in V\cap (DF\setminus\{O\})$. It follows that
$$V=(V\cap (DF\setminus\{O\})\cup(U\cap S),$$
where $V\cap (DF\setminus\{O\})\in \tau$ (hence $V\cap (DF\setminus\{O\}) = U^\prime \cap DF$ with
$U^\prime\subseteq \mathbb{A}^2_{\mathbb{R}}$ open) and $U\subseteq \mathbb{A}^2_{\mathbb{R}}$ open such that
$O\in U$. The proof of (iv) is achieved.

For property (vi), we use the fact that $S\cap W=\{O\}$ and then for an open neighborhood
$U \subset \mathbb{A}^2_{\mathbb{R}}$ of $O$ w.r.t. the real topology of $\mathbb{A}^n_{\mathbb{R}}$,
$U\cap S$ (resp. $U\cap W$) is an open neighborhood of $O\in DF$ w.r.t. $\tau_S$ (resp. $\tau_W$)
and $(U\cap S)\cap W=\{O\}$ (resp. $(U\cap W)\cap S=\{O\}$). Hence $\{O\}$ is open in $W$ (resp. in $S$)
w.r.t. $\tau_S$ (resp. $\tau_W$) and so it is a connected component of $W$ (resp. $S$),
because $\{O\}$ is also closed in $W$ (resp. $S$) w.r.t the separated topology $\tau_S$ (resp. $\tau_W$).
The connection of $p(1,\infty)$, $p(-\infty,-1)$, $S = p(-1,1)$ (resp. $p^\prime(1,\infty)$, $p^\prime(-\infty,-1)$,
$W=p^\prime(-1,1)$) w.r.t. $\tau_S$ (resp. $\tau_W$) is clear because $p$ (resp. $p^\prime$) is a homeomorphism.

{\bf Comment} (ii) in conjunction to (iii), as well as (iv), determine completely the topology $\tau_S$ (resp. $\tau_W$)
by means of the real topology of the ambient space $\mathbb{A}^2_{\mathbb{R}}$, $DF$ and its branch $S$ (resp. $W$).

\subsection{Some differential structures on affine Descartes Folium}

\hspace{0.5cm} On the topological space $(DF,\tau_S)$ (resp. $(DF,\tau_W)$) we can introduce a structure
${\cal A}_S$ (resp ${\cal A}_W$)
of smooth {\it differential manifold} by means of the simple atlas
$$\{(DF,\pi)\},\,\,(\hbox{resp.}\, \{(DF,\pi^{\prime})\}$$
having only one chart, where
$$\pi: DF \stackrel{\sim}\longrightarrow \mathbb{R}\setminus\{-1\}\subset \mathbb{R}$$
$$(\hbox{resp.}\,\,\pi^\prime: DF \stackrel{\sim}\longrightarrow \mathbb{R}\setminus\{-1\}\subset \mathbb{R})$$
is the bijective map defined above, i.e.,
$$\pi: DF \to \mathbb{R},\,\, \pi(x,y) = \left\{\begin{array}{ccc} \frac{y}{x}&if&(x,y)\not=O\\ \
0&if &(x,y)=O\end{array}\right.$$
$$\left(\hbox{resp.}\,\,\pi^{\prime}: DF \to \mathbb{R},\,\, \pi^{\prime}(x,y) = \left\{\begin{array}{ccc} \frac{x}{y}&if&(x,y)\not=O\\ \
0&if &(x,y)=O\end{array}\right.\right).$$
Recall that the inverse of the map $\pi$ (resp. $\pi^\prime$) is the map $p$ (resp. $p^\prime=\sigma p$)
and all are continuous, hence homeomorphisms. In this way,
$$
 DF\,\,\,\begin{array}{c}{\stackrel{\pi}\longrightarrow}\\ \
{\stackrel{\longleftarrow}p} \end{array}\,\,\,\mathbb{R}\setminus\{-1\}\,\,\,
\left(\hbox{resp.} DF\,\,\,\begin{array}{c}{\stackrel{\pi^\prime}\longrightarrow}\\ \
{\stackrel{\longleftarrow}{p^\prime}} \end{array}\,\,\,\mathbb{R}\setminus\{-1\} \right)
 $$
become diffeomorphisms of differentiable manifolds.

In particular, $p\alpha:\mathbb{R}\setminus \{0\}\stackrel{\sim}\longrightarrow DF$
is then also a diffeomorphism, where $DF$ is endowed with the topology $\tau_S$ and the atlas $\{(DF,\pi)\}$.
Since
$$p\alpha:(\mathbb{R}\setminus \{0\},\cdot)\stackrel{\sim}\longrightarrow (DF,\circ)$$
is a group isomorphism, it follows directly

{\bf Proposition} {\it (i) $(DF,\circ)$ is a Lie group over $\mathbb{R}$ (in particular a topological group),
where $DF$ is endowed with the topology $\tau_S$ and the differential manifold structure given by the atlas
$\{(DF,\pi)\}$.

(ii) $$p\alpha:(\mathbb{R}\setminus \{0\},\cdot)\stackrel{\sim}\longrightarrow (DF,\circ)$$
is then an isomorphism of Lie groups over $\mathbb{R}$ (in particular an isomorphism of
topological groups).}

\section{Second "Exotic" Structure}

\hspace{0.5cm} Let $\mathbb{K}=\mathbb{R}$, $DF \subset A^2_\mathbb{R}$ and $A^1_\mathbb{R}\setminus \{-1\}= \mathbb{R}\setminus \{-1\}$.
It is natural to consider also the parametrization $p^{\prime} = \sigma \circ p$. The diagram
$$
\begin{array}{cccc} DF & \left(x = \frac{3at^2}{1+t^3}, y=\frac{3at}{1+t^3}\right)&(x,y)\in DF\setminus\{(0,0)\}&(0,0)\\ \
{\wr}\Big\uparrow p^{\prime}& \Big\uparrow&\Big\downarrow&\Big\downarrow\\ \
\mathbb{R}\setminus\{-1\}&t&t=\frac{x}{y}&t=0\end{array}
$$
shows that the parametrization $p^{\prime}$ is bijective.

Since $DF\subset A^2_\mathbb{R}$, the diagram
$$
\begin{array}{cccc} DF & \left(x = \frac{3a(\tau-1)^2}{1+(\tau -1)^3}, y=\frac{3a(\tau -1)}{1+(\tau - 1)^3}\right)&(x,y)\in DF\setminus\{0\}&(0,0)\\ \
{\wr}\Big\uparrow p^{\prime}\alpha& \Big\uparrow&\Big\downarrow&\Big\downarrow\\ \
\mathbb{R}\setminus\{0\}&\tau &\tau =\frac{x}{y}+1&1\end{array}
$$
proves that $p^{\prime}\circ\alpha$ is a bijection. It follows that the group structure on $\mathbb{R}\setminus\{0\}$
transfers to $DF$:
$$(p^{\prime}\alpha)(\tau)\bot (p^{\prime}\alpha)(\tau^{\prime})\,\,\stackrel{def}= \,\,(p^{\prime}\alpha)(\tau\tau^{\prime})$$
or
$$
\left(\frac{3a(\tau -1)^2}{1+(\tau - 1)^3}, \frac{3a(\tau -1)}{1+(\tau - 1)^3}\right)\bot
\left(\frac{3a(\tau^{\prime} -1)^2}{1+(\tau^{\prime} - 1)^3}, \frac{3a(\tau^{\prime} -1)}{1+(\tau^{\prime} - 1)^3}\right)$$
$$\stackrel{def}= \,\, \left(\frac{3a(\tau\tau^{\prime} -1)^2}{1+(\tau\tau^{\prime} - 1)^3}, \frac{3a(\tau\tau^{\prime} -1)}{1+(\tau\tau^{\prime} - 1)^3}\right).$$

It appears an isomorphism of groups
$$(DF,\bot)\stackrel{p^{\prime}\alpha}{\simeq} (\mathbb{R}\setminus\{0\},\cdot)$$
$$(p^{\prime}\alpha)(\tau)\leftarrow \tau,$$
where
$$(p^{\prime}\alpha)(\tau) = \left(\frac{3a(\tau -1)^2}{1+(\tau - 1)^3},\frac{3a(\tau -1)}{1+(\tau - 1)^3}\right).$$

Here, $0=(0,0)\in DF$ is the neutral element for $(DF,\bot)$. Also, we have a canonical
isomorphism of groups over $\mathbb{R}$,
$$(DF,\circ)\stackrel{\sigma}{\stackrel{\sim}{\longrightarrow}} (DF,\bot)$$
$$(x,y)\longleftarrow (y,x),$$
i.e., $(DF,\circ)\simeq (DF,\bot)$. It is easy to see that the two group composition laws $\circ$ and $\bot$ are distinct.

{\bf Remark} As in the case of the first "exotic" structure, we can define the group structure $(DF,\bot)$,
the group isomorphism $p^\prime \alpha$ as the previous isomorphism $\sigma$ in the more general situation when the base field
$\mathbb{K}$ has $\hbox{char.}\mathbb{K}\not= 3$ and the property that the equation $\lambda^3+1=0$ has in $\mathbb{K}$
only the root $\lambda = -1$.

Let us return to the particular case when $\mathbb{K}=\mathbb{R}$.

As in the case of the first "exotic" structure, the pair  $(DF,\bot)$ is not a topological group
if $DF$ is considered with the topology $\tau$ induced by the real topology of $\mathbb{R}^2$.

By following the idea used for the first "exotic" structure, let us consider on $DF$
the topology $\tau_W$ and the topological space $(DF,\tau_W)$, the smooth differential manifold structure ${\cal A}_W$
given by the atlas $\{(DF,\pi^\prime)\}$, all being foregoing defined (see $\S$4.2 and $\S$ 4.3). Then

{\bf Proposition} {\it $(DF,\bot)$ is a Lie group over $\mathbb{R}$ (in particular, a topological group),
where $DF$ is endowed with the topology $\tau_W$ and the atlas $\{(DF,\pi^\prime)\}$.

(ii) $$p^\prime\alpha:(\mathbb{R}\setminus \{0\},\cdot)\stackrel{\sim}\longrightarrow (DF,\bot)$$
and
$$\sigma:(DF, \circ)\stackrel{\sim}\longrightarrow (DF,\bot)$$
are isomorphisms of Lie groups over $\mathbb{R}$ (in particular, isomorphisms of topological groups),
where $(DF,\circ)$ is endowed with the Lie group structure over $\mathbb{R}$, defined previously.}

It is obvious that we have a commutative diagram
$$
\begin{array}{ccccc}
(DF,\circ) & &\underset{\sim}{\overset{\sigma}\longrightarrow}& & (DF,\bot) \\
&{\wr}{\nwarrow} p\alpha &  {}  & {\wr} {\nearrow} p^{\prime}\alpha& \\  {} && (\mathbb{R}\setminus\{0\},\cdot) &&  \end{array}
$$
of isomorphisms of Lie groups over $\mathbb{R}$, as they have been defined above.

\section{First Isomorphism of \\Algebraic Lie Groups}

\hspace{0.5cm} Let $\mathbb{K}$ be a field with $\hbox{char.}\mathbb{K}\not = 3$,
$\overline{DF}\setminus\{O\} \subset \mathbb{P}^2_\mathbb{K}$, with $O=(0,0,1)$,
 $P^1_\mathbb{K}\setminus\{0,\infty\}= \mathbb{K} \setminus\{0\}$ and the infinity point $(-3a,3a,0)=(1,-1,0)$
of $\overline{DF}$. The diagram
$$
\begin{array}{ccc} \overline{DF}\setminus\{O\} & \left(x=3at, y = 3at^2, z=1+t^3\right)&(x,y,z)\in \overline{DF}\setminus\{(0,0,0)\}\\ \
{\wr}\Big\uparrow \bar{p}& \Big\uparrow&\Big\downarrow\\ \
\mathbb{K}\setminus\{0\}&t&t=\frac{y}{x}\end{array}
$$
$$
\begin{array}{c} (1,-1,0)\\ \
\Big\uparrow\\ \
t=-1\end{array}
$$
shows that the parametrization $\bar p: \mathbb{K}\setminus\{0\}\to \overline{DF}\setminus\{O\}$, $\bar p(t)=(3at,3at^2, 1+t^3)$, is bijective.
The multiplicative group on $\mathbb{K}\setminus\{0\}$ is transported on $\overline{DF}\setminus\{O\}$.
This is realized by the definition
$$\bar p(t)\cdot \bar p(t^{\prime})\,\,\overset{def}{=}\,\,\bar p(t t^{\prime})$$
or
$$(3at, 3at^2, 1+t^3)\cdot(3at^{\prime},3a{t^{\prime}}^2, 1+{t^{\prime}}^3)\,\,\overset{def}{=}\,\,(3at t^{\prime},3a({t t^{\prime}})^2, 1+({t t^{\prime}})^3),$$
for each $t, t^\prime \in \mathbb{K}\setminus\{0\}$.
Consequently, the pair $(\overline{DF}\setminus\{O\},\cdot)$ is a group and even an algebraic Lie group (see [1]) over $\mathbb{K}$,
if $\mathbb{K}$ is algebraically closed, because $(\mathbb{K}\setminus\{0\},\cdot)$ is so and $\bar p$
and ${\bar p}\,^{-1}$ are algebraic maps (morphisms).

Then, in the diagram
$$(\overline{DF}\setminus\{O\},\cdot)\stackrel{\bar p}{\simeq} (\mathbb{K}\setminus\{0\},\cdot)$$
$$\bar p(t)\longleftarrow t,$$
$\bar p$ is an isomorphism of groups. If $\mathbb{K}$ is algebraically closed, $\bar p$
is just an isomorphism of algebraic Lie groups over $\mathbb{K}$.
Finally,
$$\bar p(1) = (3a,3a,2)\in \overline{DF},$$
is the neutral element for the group
$$(\overline{DF}\setminus\{O\}, \cdot).$$

{\bf Remarks} (1) If moreover $\hbox{char.}\mathbb{K} \not = 2$, then
$$\bar p(1) = (\frac{3a}{2},\frac{3a}{2},1)\in \overline{DF},$$
equivalent to
$$\bar p(1) = (\frac{3a}{2},\frac{3a}{2})\in {DF}\subset \mathbb{A}^2_{\mathbb{K}}.$$
Let us underline that in the case $\hbox{char.}\mathbb{K} \not = 2$, the point
$\bar p(1) = (\frac{3a}{2},\frac{3a}{2})$ is just the intersection point, different of $O=(0,0)$, of $DF$
with the first bisector $x-y=0$ of $\mathbb{A}^2_{\mathbb{K}}$. In the sequel we shall denote
$$V = \bar p(1)= (3a,3a,2) \in \overline{DF}$$
(the "vertex" of $\overline{DF}$ or $DF$).

(2) If $\hbox{char.}\mathbb{K}=2$, then
$$\bar p(1) =(3a,3a,2)= (3a,3a,0) = (1,1,0)=(1,-1,0),$$
i.e., $V=I$.

Now, returning to the situation when $\hbox{char.}\mathbb{K}\not = 3$,
let $\bar p(t), \bar p(t^{\prime}) \in \overline{DF}\setminus\{O\}$, with $t, t^{\prime}\in \mathbb{K}\setminus\{0\}$.
Since $\bar p$ is a bijection, $\bar p(t), \bar p(t^{\prime})$ are any points on the curve $\overline{DF}\setminus\{O\}$.
We have
$$\bar p(t^{\prime})= [\bar p(t)]^{-1}\Longleftrightarrow \bar p(t^{\prime})\bar p(t)=\bar p(1),$$
where $\bar p$ is a morphism of groups,
$[\bar p(t)]^{-1}$ is the inverse of $\bar p(t)$ and $\bar p(1)$ is the neutral element of the group $(\overline{DF}\setminus\{O\}, \cdot)$.
The foregoing equivalence is continued by the following ones
$$\Longleftrightarrow \bar p(t^{\prime}t)=\bar p(1) \underset{\hbox{bij.}}{\stackrel{\bar p}{\Longleftrightarrow}}t^{\prime}t=1\Longleftrightarrow t^{\prime}=\frac{1}{t}\Longrightarrow [\bar p(t)]^{-1}=\bar p\left(\frac{1}{t}\right), \forall t \in \mathbb{K}\setminus\{0\}.$$

\subsection{Geometrical interpretation}

\hspace{0.5cm} Let $\mathbb{K}$ be an arbitrary field and
$$\begin{array}{ccc} A^2_\mathbb{K} &\stackrel{\sigma}{\stackrel{\sim}{\longrightarrow}} &A^2_\mathbb{K}\\ \
(x,y)&\longrightarrow &(y,x)\end{array}$$
the symmetry of $A^2_\mathbb{K}$ w.r.t. the first bisector $x-y=0$.
The application $\sigma$ is bijective and $\sigma(DF)=DF$.
We extend $\sigma$ to the bijective map
$$\begin{array}{ccc} P^2_\mathbb{K} &\stackrel{\bar\sigma}{\stackrel{\sim}{\longrightarrow}} &P^2_\mathbb{K}\\ \
(x,y,z)&\longrightarrow &(y,x,z).\end{array}$$
We have $\bar\sigma(\overline{DF})= \overline{DF}$.

Let us point out that for each field $\mathbb{K}$, with $\hbox{char.}\mathbb{K}\not = 3$, the following diagram
$$\begin{array}{ccc}\overline{DF}\setminus\{O\}& \underset{\sim}{\overset{\bar\sigma}\longrightarrow}&\overline{DF}\setminus\{O\}\\ \
{\bar p}\Big\uparrow  \wr& &{\bar p}\Big\uparrow  \wr\\ \ \mathbb{K}\setminus\{0\}& \overset{\sim}\longrightarrow &\mathbb{K}\setminus\{0\}\\ \
t&\longrightarrow  &\frac{1}{t}\end{array}$$
is commutative, i.e., $\bar \sigma \bar p(t) = \bar p\left(\frac{1}{t}\right)$, for each $t\in \mathbb{K}\setminus\{0\}$.

Indeed, for $t\in \mathbb{K}\setminus\{0\}$, we have
$$\bar\sigma \bar p(t) = \bar\sigma(3at, 3at^2, 1+t^3)$$
$$= (3at^2,3at, 1+t^3) = \left(\frac{3a}{t}, \frac{3a}{t^2}, 1 + \frac{1}{t^3}\right)=\bar p \left(\frac{1}{t}\right).$$
In particular, for each $t\in \mathbb{K}\setminus\{0\}$, we have
$\bar\sigma \bar p(t)=\left[\bar p(t)\right]^{-1}$.
If $\bar p(t)=(x,y,z)\in \overline{DF}\setminus\{O\}$, then
$$(x,y,z)^{-1}= \left[\bar p(t)\right]^{-1}= \bar\sigma \bar p(t)=\bar \sigma(x,y,z)= (y,x,z)$$
(i.e., another writing of the inverse w.r.t the composition law $\cdot$).

In the particular case when $z=1$, we have that for each point $(x,y)\in DF\setminus\{O\}\subset \mathbb{A}^2_{\mathbb{K}}$,
its symmetric/opposite w.r.t. the composition law $\cdot$ is $(x,y)^{-1}=(y,x)\in DF\setminus\{O\}\subset \mathbb{A}^2_{\mathbb{K}}$,
i.e., the symmetric of the point $(x,y)$ w.r.t. the first bisector $x-y=0$ of $\mathbb{A}^2_{\mathbb{K}}$.

If $t=-1$, then $p(-1)=(-3a,3a,0)= (1,-1,0)$ is the point at infinity $I$ of $\overline{DF}\setminus\{O\}$ and then
$$I^{-1}=[p(-1)]^{-1}= [(1,-1,0)]^{-1}= (-1,1,0)= (1,-1,0) = p(-1)=I.$$

\section{Second Isomorphism of \\Algebraic Lie Groups}

\hspace{0.5cm} Let $\mathbb{K}$ be a field with $\hbox{char.}\mathbb{K}\not=3$ and $P^1_\mathbb{K}\setminus\{0,\infty\}=\mathbb{K}\setminus \{0\}$.
As we showed, the diagram
$$\begin{array}{ccc} \mathbb{P}^2_\mathbb{K} &\stackrel{\bar\sigma}{\stackrel{\sim}{\longrightarrow}} &\mathbb{P}^2_\mathbb{K}\\ \
(x,y,z)&\longrightarrow &(y,x,z)\end{array}$$
implies
$$\bar\sigma(\overline{DF})= \overline{DF}.$$
Let us consider the parametrization $\bar{\bar{p}}$ and the diagram
$$
\begin{array}{cccc} \overline{DF}\setminus\{O\} & (x=3at^2, y = 3at, z=1+t^3)&(x,y,z)&(1,-1,0)\\ \
{\wr}\Big\uparrow \bar{\bar{p}}& \Big\uparrow&\Big\downarrow&\Big\uparrow\\ \
\mathbb{K}\setminus\{0\}&t&t=\frac{x}{y}&t=-1\,.\end{array}
$$
The parametrization
$\bar{\bar p}: \mathbb{K}\setminus\{0\}\to \overline{DF}\setminus\{O\}$, $\bar{\bar p}(t)=(3at^2,3at, 1+t^3)$, is bijective.
The multiplicative group structure group on $\mathbb{K}\setminus\{0\}$ can be transported on $\overline{DF}\setminus\{O\}$.
This is realized by the definition
$$\bar{\bar p}(t)\circ \bar{\bar p}(t^{\prime})\,\,\overset{def}{=}\,\,\bar{\bar p}(t \cdot t^{\prime})$$
or
$$(3at^2, 3at, 1+t^3)\cdot(3a{t^{\prime}}^2,3a{t^{\prime}}, 1+{t^{\prime}}^3)\,\,\overset{def}{=}\,\,(3a({t t^{\prime}})^2,3a({t t^{\prime}}), 1+({t t^{\prime}})^3),$$
for each $t, t^\prime\in \mathbb{K}\setminus\{0\}$. Consequently,
the pair $(\overline{DF}\setminus\{O\},\circ)$ is a group and even an algebraic Lie group (see [1]) over $\mathbb{K}$
if $\mathbb{K}$ is algebraically closed,
because $(\mathbb{K}\setminus\{0\},\cdot)$ is so, and $\bar{\bar p}$ and ${\bar{\bar p}}\,^{-1}$ are algebraic maps.
Therefore in the diagram
$$(\overline{DF}\setminus\{O\},\circ)\stackrel{{\bar{\bar p}}}{\simeq} (\mathbb{K}\setminus\{0\},\cdot)$$
$$\bar{\bar p}(t)\longleftarrow t$$
$\bar{\bar p}$ is an isomorphism of groups and in the situation when $\mathbb{K}$ is algebraically closed,
$\bar{\bar p}$ is just an isomorphism of algebraic Lie groups over $\mathbb{K}$. Finally,
$$\bar{\bar p}(1) = (3a,3a,2)= V\in \overline{DF},$$
(which is equivalent to
$$\bar{\bar p}(1)=V=\left(\frac{3a}{2},\frac{3a}{2}\right)\in {DF}\subset \overline{DF}$$
if $\hbox{char.}\,\mathbb{K} \not = 2,3$).
In this way, the point $V= \bar{\bar p}(1)\in \overline{DF}$ is the neutral element for the group
$$(\overline{DF}\setminus\{O\}, \circ).$$
Therefore the neutral elements of $(\overline{DF}\setminus\{O\}, \cdot)$ and $(\overline{DF}\setminus\{O\}, \circ)$ coincides.
As in Section 6, we have the inversion formula with respect to the composition law $\circ$:
$$\left[\bar{\bar p}(t)\right]^{-1}(t) = \bar \sigma ({\bar{\bar p}}(t)),\forall t \in \mathbb{K}\setminus\{0\}$$
and, for each $(x,y,z)\in \overline{DF}\setminus\{O\}$, we have $(x,y,z)^{-1}= (y,x,z).$ This means that, for any point $\bar{\bar p}(t) \in (\overline{DF}\setminus\{O\},\circ)$,
the inverses with respect to each of operations of groups $\circ$ and $\cdot$ coincide.

We have the following commutative diagram
$$\begin{array}{ccccc}(\overline{DF}\setminus\{O\},\cdot)& &\underset{\sim}{\overset{\bar\sigma}\longrightarrow}&&(\overline{DF}\setminus\{O\},\circ)\\ \
&{\wr}{\nwarrow} \bar{p} & &{\wr}{\nearrow} \bar{\bar p}& \\ \ &&(\mathbb{K}\setminus\{0\},\cdot)&& \end{array}$$
Since $\bar p$ and $\bar{\bar p}$ are both isomorphisms of (algebraic Lie) groups (over $\mathbb{K}$),
the function $\bar \sigma$ is also such an isomorphism. According to the previous remark on the inversion formula,
for each $A\in \overline{DF}\setminus\{O\}$, we have $\bar \sigma(A)= A^{-1}$.

{\bf Proposition}. The groups $(\overline{DF}\setminus\{O\},\cdot)$ and $(\overline{DF}\setminus\{O\},\circ)$
coincide, i.e. $\cdot = \circ$.

{\bf Proof} For all $A,B \in \overline{DF}\setminus\{O\}$, we have
$$(A^{-1}\cdot B^{-1})^{-1}= \bar{\sigma}(A^{-1}\cdot B^{-1}),$$
which splits as
$$(A^{-1}\cdot B^{-1})^{-1}=(A^{-1})^{-1}\cdot (B^{-1})^{-1}=A\cdot B$$
$$ \bar{\sigma}(A^{-1}\cdot B^{-1})= \bar{\sigma}(A^{-1})\circ \bar{\sigma}(B^{-1})= (A^{-1})^{-1}\circ (B^{-1})^{-1}= A\circ B.$$
Here we used the fact that the inverse $A^{-1}$ of a point $A$ does not depend on the operation
$\cdot$ or $\circ$ with respect to which we consider it.

{\bf Lemma} {\it Let $\mathbb{K}$ be a field and
$$\overline{DF}: x^3 +y^3 -3axyz=0\, \subset \mathbb{P}^2_{\mathbb{K}},\, a\in \mathbb{K}\setminus\{0\}.$$
Let us consider a straight line $\bar d\subset \mathbb{P}^2_{\mathbb{K}}$ which cuts $\overline{DF}$ in the points
$$P_1(x_1,y_1,z_1), P_2(x_2,y_2,z_2), P_3(x_3,y_3,z_3)\in \mathbb{P}_{\mathbb{K}}^2$$
(counted with multiplicity). Then
$$x_1x_2x_3 + y_1y_2y_3=0.$$
In particular, if $x_1x_2x_3\not=0$, then $t_1t_2t_3=-1$, where
$t_i=\frac{y_i}{x_i}$ is the slope of the affine straight line $\overline{OP}_i\subset \mathbb{P}_{\mathbb{K}}^2$
on $\mathbb{A}^2_{\mathbb{K}}\subseteq \mathbb{P}^2_{\mathbb{K}}$.}

{\bf Proof} Let us suppose that
$$\bar d: mx+ny + pz=0 \subset \mathbb{P}^2_{\mathbb{K}}.$$
It is obvious that $O(0,0,1)\in \bar d \Longleftrightarrow p=0$.

If $O\in \bar d$, then we may assume that $P_1=O$. Then $x_1=y_1=0$ and the relation from
the Lemma is fulfilled.

If $O\notin \bar d$, then $p\not=0$ and we may assume that $p=-1$, i.e.,
$$\bar d: mx+ny =z\subset \mathbb{P}^2_{\mathbb{K}}.$$
In this case $P_1, P_2, P_3 \not = O$ and $x_i\not = 0$, for all $i=1,2,3$.
Since $\{P_1,P_2,P_3\}= \overline{DF}\cap \bar d$, the pairs $(x_i,y_i), i=1,2,3$ verify the relation
$$x^3 + y^3 -3axy(mx + ny)=0$$
and $t_i=\frac{y_i}{x_i}, i=1,2,3$, are the roots of the equation
$$t^3 -3ant^2 -3amt+1=0.$$
Consequently $t_1t_2t_3 = -1$, i.e., $x_1x_2x_3 + y_1y_2y_3=0$.

\vspace{0.2cm}
{\bf Lemma} {\it Let $\mathbb{K}$ be a field with $\hbox{char.}\,\mathbb{K}\not = 3$ and
$$\overline{DF}: x^3 +y^3 -3axyz=0\, \subset \mathbb{P}^2_{\mathbb{K}},\,a\in \mathbb{K}\setminus\{0\}$$
$$P_1=\bar p(t_1), P_2=\bar p(t_2), P_3=\bar p(t_3)\in \overline{DF}\setminus\{0\},$$
with $t_1,t_2,t_3\in \mathbb{K}\setminus\{0\}$. Then $P_1, P_2, P_3$ are the intersections
(counted with multiplicity) of a straight line
$\bar d \subset \mathbb{P}^2_{\mathbb{K}}$ with $\overline{DF}\setminus\{0\}$ if and only if $t_1t_2t_3=-1$.}

{\bf Proof} ($\Longrightarrow$) Let $P(x,y,z))\in \overline{DF}\setminus\{O\}$, with $P=\bar p(t)$, where
$t\in \mathbb{K}\setminus\{0\}$. Then $t= {\bar p}^{-1}(P)=\frac{y}{x}$, according to the definition of ${\bar p}^{-1}$.
If
$$P_1=\bar p(t_1), P_2=\bar p(t_2), P_3= \bar p(t_3)\in \overline{DF}\setminus\{O\}$$
are the intersection points (counted with multiplicity) of $\overline{DF}\setminus\{O\}$ with the straight line
$\bar d \subset \mathbb{P}^2_{\mathbb{K}}$, then it is obvious that $t_1t_2t_3=-1$, according to the previous Lemma.

($\Longleftarrow$) To prove the converse assertion, let us consider the straight line
$\bar d\subset \mathbb{P}^2_{\mathbb{R}}$ determined by $P_1,P_2$ (which is particularly the tangent line to
$\overline{DF}$ at $P_1=P_2$, if these points coincide). Let $P_3^\prime$ be the third point
of intersection of $\bar d$ with $\overline{DF}$ (counted with multiplicity). Suppose $P_3^\prime=\bar p(t_3^\prime)$,
with $t_3^\prime\in \mathbb{K}\setminus\{0\}$. Then, according to the first part of this proof, we have $t_1t_2t_3^\prime=-1$.
Because $t_1t_2t_3=-1$ and $t_1\not= 0, t_2\not=0$, it follows $t_3^\prime=t_3$ and $P_3^\prime=P_3$. Consequently
$P_1, P_2, P_3$ are the intersection points of $\bar d$ with $\overline{DF}$ (counted with the multiplicity).

We can formulate an equivalent form of the previous Lemma involving the group structure
$\cdot$ on $\overline{DF}\setminus\{0\}$:
\vspace{0.2cm}

{\bf Theorem} {\it Let $\mathbb{K}$ be a field with $\hbox{char.}\,\mathbb{K} \not = 3$ and
$$\overline{DF}: x^3 +y^3 -3axyz=0\, \subset \mathbb{P}^2_{\mathbb{K}},\,a\in \mathbb{K}\setminus\{0\},$$
$$I=(1,-1,0)\in \overline{DF}\setminus\{O\},\,\hbox{the point at infinity of}\,\,\,\overline{DF},$$
$$P_1,P_2,P_3\in \overline{DF}\setminus\{O\}.$$
Then $P_1,P_2,P_3$ are the intersection points of a straight line $\bar d\subset \mathbb{P}^2_{\mathbb{K}}$
with $\overline{DF}\setminus\{O\}$ if and only if $P_1\cdot P_2\cdot P_3=I$.}

{\bf Proof} Suppose $P_i=\bar p(t_i)$, with $t_i\in \mathbb{K}\setminus\{0\},\, i=1,2,3$. Recall that
$I=\bar p(-1)$. We have $t_1t_2t_3=-1$ if and only if $\bar p(t_1t_2t_3)=\bar p(-1)$ or if and only if
$\bar p(t_1)\bar p(t_2)\bar p(t_3)=I$, which achieved the proof.

{\bf Definition} {\it Let $\mathbb{K}$ be a field with $\hbox{char.}\,\mathbb{K}\not = 3$ and $P=p(t)\in \overline{DF}\setminus\{O\}$
with $t\in \mathbb{K}\setminus\{0\}$. We define $P^{\bot}\in \overline{DF}\setminus\{O\}$ by the relation $P^\bot=\bar p\left(-\frac{1}{t}\right)$}.

{\bf Remarks} 1) Suppose $\mathbb{K}=\mathbb{R}$. Then the slope of the affine trace of the straight line
$\overline{OP}\subset \mathbb{P}^2_{\mathbb{R}}$, resp. $\overline{OP^\bot} \subset \mathbb{P}^2_{\mathbb{R}}$,
on $\mathbb{A}^2_{\mathbb{R}}\subset \mathbb{P}^2_{\mathbb{R}}$, is $t$, resp. $-\frac{1}{t}$,
hence they are orthogonal vectorial lines in $\mathbb{A}^2_{\mathbb{R}} = \mathbb{R}^2$ w.r.t. the canonical Euclidean structure.
Therefore $P^\bot$ receives a geometric definition if $\mathbb{K}=\mathbb{R}$.

2) We have ${\bar p}(t)^\bot = \bar p\left(-\frac{1}{t}\right)$ for each $t \in \mathbb{K}\setminus\{0\}$. It follows that:

(i) $(P^\bot)^\bot = P$ for each $t\in \mathbb{K}\setminus\{0\}$;

(ii) if $P=I=(1,-1,0) = \bar p(-1)$, then $I^\bot = \bar p\left(-\frac{1}{-1}\right)=\bar p(-1)=V$ and consequently $V^\bot =I$;

(iii) if $P\in DF\setminus\{O\},\,P\not = V$ (equivalently $P=\bar p(t)$ with $t\in \mathbb{K}\setminus\{0,1\}$), then
$P^\bot \in DF\setminus\{O\}$ and conversely.

Now we can give a geometric definition of the composition law $\cdot$ of the group $(\overline{DF}\setminus\{0\},\cdot)$
as follows:

{\bf Theorem} {\it Let $\mathbb{K}$ be a field with $\hbox{char.}\,\mathbb{K}\not =3$ and
$$\overline{DF}: x^3+y^3-3axyz=0\subset \mathbb{P}^2_{\mathbb{K}},\, a \in \mathbb{K}\setminus\{0\}.$$
Suppose that:(i) $P_1, P_2\in \overline{DF}\setminus\{O\}$ are distinct (resp. non distinct) points;
(ii) $\overline{P_1P_2}$ is the straight line (resp. the tangent line to $\overline{DF}$)
determined by $P_1, P_2$ in $\mathbb{P}^2_{\mathbb{K}}$;
(iii) $P_3 \in \overline{DF}\setminus\{O\}$ (counted with multiplicity) is the third intersection point of
$\overline{P_1P_2}$ with $\overline{DF}\setminus\{O\}$.

Then
$$P_1\cdot P_2= P_3^\bot.$$}

{\bf Proof} Suppose $P_1=\bar p(t_1), P_2=\bar p(t_2), P_3=\bar p(t_3)$ with $t_1, t_2, t_3 \in \mathbb{R}\setminus\{0\}$.
By the previous Lemma, we have $t_1t_2t_3=-1$. Then
$$P_1\cdot P_2= \bar p(t_1)\bar p(t_2)= \bar p(t_1t_2) = \bar p\left(-\frac{1}{t_3}\right)=\bar p(t_3)^\bot = P_3^\bot.$$

As application of this definition of the group law $\cdot$, for $\mathbb{K}=\mathbb{R}$, we have the
following pure geometric property of affine Descartes Folium, not involving any group structure.

{\bf Corollary} {\it Let $\mathbb{K}=\mathbb{R}$ and
$$DF: x^3+y^3-3axy=0\subset \mathbb{A}^2_{\mathbb{R}},\, a \in \mathbb{R}\setminus\{0\}.$$
We fix the points $V =\left(\frac{3a}{2},\frac{3a}{2}\right)\in DF$, $P\in DF\setminus\{O\}, P\not = V$,
and let $Q\in DF\setminus\{O\}$ be the third intersection point of the affine straight line $VP$ with $DF$. Then $OP\bot OQ$
(perpendicular straight lines in $\mathbb{A}^2_{\mathbb{R}}=\mathbb{R}^2$ w.r.t. the canonical Euclidean structure).

Conversely, if $P,Q \in DF\setminus\{O\}$, with $P,Q\not= V$ and $OP\bot OQ$, then $P, Q, V$ are collinear points.}

{\bf Proof} Since $V$ is the neutral element of the group $(\overline{DF}\setminus\{O\}, \cdot)$, we have $V\cdot P = P$.
On the other hand, according to the Theorem above, $V\cdot P = Q^\bot$. Consequently, $P=Q^\bot$ and then $OP\bot OQ$
by the definition of $Q^\bot$.

Conversely, suppose $OP\bot OQ$ and let $R$ be the third intersection point
of the straight line $VQ$ with $DF\setminus\{0\}$. Then $OR\bot OQ$. It follows $OP=OR$, $P=R$ and so $P=R\in VQ$.

Now, in conjunction with the Corollary above, the previous Theorem about a geometric definition of $\cdot$,
can be rewritten in the following form:

{\bf Theorem} {\it Let $\mathbb{K}$ be a field with $\hbox{char.}\,\mathbb{K}\not = 3$, let $\overline{DF}$ be the
projective Descartes Folium over $\mathbb{K}$ and $P_1, P_2, P_3 \in \overline {DF}\setminus\{O\}$
be as in the previous Theorem. Denote by $Q$ the third intersection point of the line
$\overline{VP_3}\subset \mathbb{P}^2_{\mathbb{K}}$ with $\overline {DF}\setminus\{O\}$
(counted with multiplicity). Then $P_1\cdot P_2 = Q$.}

{\bf Remark} If we have in mind that
$$V=(3a,3a,2)\in \overline {DF}\setminus\{O\}\subset \mathbb{P}^2_{\mathbb{K}}$$
is the neutral element of the group $(\overline {DF}\setminus\{O\}, \cdot)$,
this last geometric definition of the composition law $\cdot$ in the group
$(\overline {DF}\setminus\{O\}, \cdot)$ recall the classic well-known geometric definition of
the group composition law on elliptic curves.

{\bf Proof of Theorem} Suppose
$$P_1=\bar p(t_1), P_2=\bar p(t_2), P_3=\bar p(t_3), Q=\bar p(t),$$
with $t_1,t_2,t_3,t \in \mathbb{K}\setminus\{0\}$. Recall that $V=\bar p(1)$.
Since $P_1,P_2,P_3$, resp. $V,P_3,Q$, are the intersection points of a straight line with
$\overline {DF}\setminus\{0\}$, by a previous Lemma we have $t_1t_2t_3=-1$, resp. $1\cdot t_3\cdot t=-1$.
It follows $t_1t_2=t$ and hence
$$P_1\cdot P_2=\bar p(t_1)\cdot\bar p(t_2)= \bar p(t_1t_2)=\bar p(t)=Q.$$

\section{A Derived Structure $(\overline {DF}\setminus\{0\},\star )$ of \\Algebraic Lie Group}

\hspace{0.5cm} Let $\mathbb{K}$ be a field with $\hbox{char.}\, \mathbb{K}\not=3$.
We have a bijective map, which is an isomorphism of algebraic varieties over $\mathbb{K}$:
$$?^\bot: \overline {DF}\setminus\{O\}\stackrel{\sim}{\longrightarrow} \overline {DF}\setminus\{O\}$$
$$P=\bar p(t)\longrightarrow  P^\bot=\bar p\left(-\frac{1}{t}\right)$$
$$Q^\bot \longleftarrow Q.$$

Then the group composition law $\cdot$ can be transported by this bijective map and
we obtain a new group composition law $\star$ on $\overline {DF}\setminus\{O\}$:
$$\forall, P, P^\prime \in \overline {DF}\setminus\{O\},\,\, P\star P^\prime\,\,\overset{def}{=}\,\, (P^\bot\cdot {P^\prime}^\bot)^\bot.$$
Equivalently,
$$\bar p(t)\star\bar p(t^\prime)\,\,\overset{def}{=}\,\, \bar p(-tt^\prime),\, \forall t, t^\prime \in \mathbb{K}\setminus\{0\}$$
or
$$\bar p(-t)\star\bar p(-t^\prime)\,\,\overset{def}{=}\,\, \bar p(-tt^\prime),\, \forall t, t^\prime \in \mathbb{K}\setminus\{0\}$$

Thus $(\overline {DF}\setminus\{O\},\star)$ is a group and if $\mathbb{K}$ is algebraically closed,
it becomes an algebraic Lie group over $\mathbb{K}$ (see [1]).
We have an isomorphism of (algebraic Lie) groups (over $\mathbb{K}$), represented schematically by
$$?^\bot: (\overline {DF}\setminus\{O\},\cdot)\stackrel{\sim}{\longrightarrow} (\overline {DF}\setminus\{O\},\star)$$
$$P=\bar p(t)\longrightarrow  P^\bot=\bar p\left(-\frac{1}{t}\right)$$
$$Q^\bot \longleftarrow Q.$$

By composition with
$$\bar p:(\mathbb{K}\setminus\{0\}, \cdot)\longrightarrow (\overline {DF}\setminus\{O\},\cdot)$$
$$t \longrightarrow \bar p(t),$$
we obtain the isomorphism of (algebraic Lie) groups (over $\mathbb{K}$), which represents also a parametrization
of the curve $\overline {DF}\setminus\{O\}$, namely,
$${\bar p}^\bot = ?^{\bot}\circ \bar p: (\mathbb{K}\setminus\{0\},\cdot)\longrightarrow (\overline {DF}\setminus\{O\},\star)$$
$$t \longrightarrow {\bar p(t)}^\bot = \bar p\left(-\frac{1}{t}\right)=\left(-\frac{3a}{t}, \frac{3a}{t^2}, 1-\frac{1}{t^3}\right)=
(-3at^2,3at,t^3-1)$$
and then
$$(-3at^2,3at,t^3-1)\star(-3a{t^\prime}^2, 3at^\prime, {t^\prime}^3-1)= (-3a(tt^\prime)^2, 3att^\prime, (tt^\prime)^3-1),$$
$$\forall t, t^\prime\in \mathbb{K}\setminus\{0\}.$$

Let us illustrate that this composition law is just that from the paper [11].
Indeed, it is easy to establish some properties of the group $(\overline {DF}\setminus\{O\},\star)$
confirming this fact:

(a) $I=(1,-1,0)\in \overline {DF}\setminus\{0\}$is the neutral element w.r.t. the composition law $\star$.
In fact $?^\bot$ is a morphism of groups, $V=\bar p(1)$ is the neutral element of the group $(\overline {DF}\setminus\{O\},\cdot)$
and $I=V^\bot$.

(b) For each $P\in \overline {DF}\setminus\{O\}$, the symmetric/opposite element $P^{-1}$
w.r.t. the group law $\star$ is the symmetric of $P$ w.r.t. the first bisector of $\mathbb{A}^2_{\mathbb{K}}$
(i.e., the symmetric/opposite elements of $P$ w.r.t. $\star$ and $\cdot$ coincide). Indeed, if
$P=\bar p(t)\in \overline {DF}\setminus\{O\}$ and if we consider $P^\prime=\bar p\left(\frac{1}{t}\right)$,
the symmetric of $P$ w.r.t. the first bisector of $\mathbb{A}^2_{\mathbb{K}}$, then we have
$$P\star P^\prime = \bar p(t)\star \bar p\left(\frac{1}{t}\right)= \bar p\left(-t\cdot\frac{1}{t}\right)=\bar p(-1)=I.$$

(c) {\bf Geometric definition of the group composition law $\star$} (see [9]) Let $P_1,P_2\in \overline {DF}\setminus\{O\}$
be distinct points (resp. not distinct points), $\bar d$ the straight lines $\overline{P_1P_2}$
(resp. the tangent line to $\overline {DF}\setminus\{O\}$ at $P_1=P_2$) and $P_3$ the third
intersection point of $\bar d$ with $\overline {DF}\setminus\{O\}$ (counted with multiplicity).
Then $P_1\star P_2$ is the symmetric of $P_3$ w.r.t. the first bisector of $\mathbb{A}^2_{\mathbb{K}}$.
In fact, if
$$P_1=\bar p(t_1), P_2=\bar p(t_2),P_3=\bar p(t_3),$$
then we have $t_1t_2t_3=-1$ and consequently
$$P_1\star P_2= \bar p(t_1)\star \bar p(t_2) \underset{\star}{\stackrel{def}{=}}\bar p(-t_1t_2)=\bar p\left(\frac{1}{t_3}\right)$$
i.e., $P_1\star P_2$ is the symmetric of $P_3=\bar p(t_3)$ w.r.t. the first bisector.

(d) Let $P_1,P_2,P_3\in \overline {DF}\setminus\{O\}$. We have $P_1\star P_2\star P_3= P_1\cdot P_2\cdot P_3$.
Then $P_1,P_2,P_3$ are the intersection points of a straight line $\bar d \subset \mathbb{P}^2_{\mathbb{K}}$
if and only if $P_1\star P_2\star P_3 = I$. Indeed, suppose $P_i=\bar p(t_i)$ with $t_i\in \mathbb{K}\setminus\{0\},\, i=1,2,3$.
Then
$$P_1\star P_2\star P_3 = \bar p(t_1)\star \bar p(t_2)\star \bar p(t_3) \underset{\star}{\stackrel{def}{=}}\bar p(-t_1t_2)\star \bar p(t_3)$$
$$\underset{\star}{\stackrel{def}{=}}\bar p(t_1t_2t_3) \underset{\cdot}{=} \bar p(t_1)\cdot \bar p(t_2)\cdot \bar p(t_3)=P_1\cdot P_2\cdot P_3.$$
The stated equivalence is a direct consequence for the similar previous property stated for the composition law $\cdot$.

(e) Let
$$P_1= {\bar p}^\bot(t_1),P_2= {\bar p}^\bot(t_2),P_3= {\bar p}^\bot(t_3)\in \overline {DF}\setminus\{O\}$$
with $t_1,t_2,t_3 \in \mathbb{K}\setminus\{0\}$. Then $P_1,P_2,P_3$ are the intersection points of a straight
line $\bar d \subset \mathbb{P}^2_{\mathbb{R}}$ if and only if $t_1t_2t_3=1$. According to (d),
to prove this statement, it is enough
to verify that $P_1\star P_2\star P_3=I$ if and only if $t_1t_2t_3=1$. In fact,
$$P_1\star P_2\star P_3= {\bar p}^\bot(t_1)\star {\bar p}^\bot(t_2)\star {\bar p}^\bot(t_3)
\underset{morph.}{\stackrel{{\bar p}^\bot}{==}}{\bar p}^\bot(t_1t_2t_3)$$
and so $P_1\star P_2\star P_3=I$ if and only if ${\bar p}^\bot(t_1t_2t_3)={\bar p}^\bot(1)$, if and only if
$t_1t_2t_3=1$.

(f) {\bf Other relations between the composition laws $\cdot$ and $\star$}

Recall the notations $V=(3a,3a,2)$, $I=(1,-1,0)$, which are the neutral elements w.r.t.
the composition laws $\cdot$, resp. $\star$. We have, for each $P_1,P_2 \in \overline {DF}\setminus\{O\}$:
$$P_1\star P_2 = P_1\cdot P_2\cdot I,\,\,\, P_1\cdot P_2 = P_1\star P_2\star V.$$
In fact, these relations are particular cases of the relation $P_1\star P_2\star P_3= P_1\cdot P_2\cdot P_3$,
when $P_3=I \,\hbox{or}\, V$.
More general, for each $n\geq 1$, we find
$$P_1\star P_2\star\ldots \star P_{2n}= P_1\cdot P_2\cdot \ldots \cdot P_{2n}\cdot I$$
$$P_1\cdot P_2\cdot\ldots \cdot P_{2n}= P_1\star P_2\star \ldots \star P_{2n}\star V$$
$$P_1\star P_2\star\ldots \star P_{2n+1}= P_1\cdot P_2\cdot \ldots \cdot P_{2n+1}.$$
For the first relation, we proceed by induction on $n\geq 1$: if $n\geq 2$, we have
$$P_1\star P_2\star\ldots \star P_{2n}= (P_1\star P_2\star\ldots \star P_{2n-2})\star(P_{2n-1}\star P_{2n})$$
$$=(P_1\cdot P_2\cdot\ldots \cdot P_{2n-2}\cdot I)\star(P_{2n-1}\cdot P_{2n}\cdot I)$$
$$=P_1\cdot P_2\cdot\ldots \cdot P_{2n-2}\cdot I\cdot P_{2n-1}\cdot P_{2n}\cdot I\cdot I$$
$$= P_1\cdot P_2\cdot \ldots \cdot P_{2n}\cdot I,$$
because $I\cdot I\cdot I=I^3=I$.

In a similar way can be established the second relation.

In the case of the third relation, for $n\geq 1$, we have
$$P_1\star P_2\star\ldots \star P_{2n+1}= (P_1\star P_2\star\ldots \star P_{2n}) \star P_{2n+1}$$
$$=(P_1\cdot P_2\cdot\ldots \star P_{2n}\cdot I) \star P_{2n+1}= P_1\cdot P_2\cdot \ldots \cdot P_{2n}\cdot I\cdot P_{2n+1}\cdot I
=P_1\cdot P_2\cdot \ldots \cdot P_{2n+1},$$
because $I\cdot I=I^2=V$.

(g) {\bf Expressions of the map $?^\bot$ w.r.t. the composition laws $\cdot$ and $\star$}

For each $P\in \overline {DF}\setminus\{O\}$, we have $P^\bot = P^{-1}\cdot I = P^{-1}\star V$.
To check this, let us suppose that $P=\bar p(t)$. Then
$$P^\bot =\bar p(t)^\bot = \bar p\left(-\frac{1}{t}\right)= \bar p\left(\frac{1}{t}\right)\cdot \bar p(-1)=P^{-1}\cdot I$$
and
$$P^\bot = \bar p\left(-\frac{1}{t}\right)= \bar p\left(\frac{1}{t}\right)\star \bar p(1)=P^{-1}\star V.$$

\section{Third "Exotic" Structure}

\hspace{0.5cm} Let $\mathbb{K}$ be a field with $\hbox{char.}\,\mathbb{K}\not=3$.
Then $\overline{DF}\subset \mathbb{P}^2_{\mathbb{K}}$ has one or three points at infinite. The diagram
$$
\begin{array}{ccc} \overline{DF} & (x=3at, y = 3at^2, z=1+t^3)&(x,y,z)\not=O\\ \
\bar p\Big\uparrow \wr&\Big \uparrow&\Big\downarrow\\ \
\mathbb{A}^1_\mathbb{K}=\mathbb{K}&t&t=\frac{y}{x}\end{array}
$$
$$
\begin{array}{ccc} (0,0,1)&(1,-1,0)&(1,\epsilon_i,0)\\ \
\Big\downarrow \Big\uparrow&\Big\uparrow&\Big\uparrow\\ \
t=0&t=-1&t=\epsilon_i\end{array}
$$
shows that the parametrization $\bar p$ is bijective. The function $\bar p$ transports the
additive group structure from $\mathbb{K}$ to $\overline{DF}$ and we obtain the group $(\overline{DF},+)$ defined by
$$\bar p(t) +  \bar p(t^{\prime}) \,\,\overset{def}{=}\,\, \bar p(t+t^{\prime})$$
or, equivalently, by
$$(3at,3at^2,1+t^3)+ (3at^{\prime},3a{t^{\prime}}^2,1+{t^{\prime}}^3)\,\,\overset{def}{=}\,\, (3a(t+t^{\prime}),3a(t+{t^{\prime}})^2,1+(t+{t^{\prime}})^3).$$
It appears an isomorphism of groups
$$(\overline{DF}, +)\stackrel{\bar p}{\stackrel{\sim}{\longleftarrow}}(\mathbb{K}, +)$$
$$\bar p(t)\longleftarrow t$$

The point $\bar p(0) = (0,0,1)=O$ is a neutral element for the group $(\overline{DF}, +)$.
We have
$$\forall t \in \mathbb{K},\,\, -\bar p(t)= \bar p(-t).$$
In this way, $-\bar p(-1)= \bar p(1)$ implies the identification between $-(1,-1,0)$ (the opposite of the point $I$ at infinite)
and $V=(3a,3a,2)\in\overline{DF}$ (the "vertex" of $\overline{DF}\subset \mathbb{P}^2_{\mathbb{K}}$).

If $\mathbb{K}=\mathbb{R}, \mathbb{C}$, as in the case of the previous "exotic" structures, the pair $(\overline{DF},+)$ is not a
topological group, where $\overline{DF}$ is considered with topology induced by the real, resp. complex,
topology of $\mathbb{P}^2_{\mathbb{R}}$, resp. $\mathbb{P}^2_{\mathbb{C}}$.
We will present in the sequel section a "correction" of this situation.

\section{Fourth "Exotic" Structure}

\hspace{0.5cm} Let $\mathbb{K}$ be a field with $\hbox{char.}\,\mathbb{K}\not = 3$. Then $\overline{DF}\subset \mathbb{P}^2_{\mathbb{K}}$
has one or three points at infinite.
Denote $\bar \sigma\circ \bar p=\bar{\bar p}$. Recall that
$$\mathbb{P}^2_\mathbb{K} \stackrel{\bar \sigma}{\stackrel{\sim}{\longrightarrow}}\mathbb{P}_\mathbb{K}^2,\,\,\bar\sigma (x,y,z)=(y,x,z)$$
$$\bar\sigma(\overline{DF}) = \overline{DF}.$$

Because the roots $\epsilon_1, \epsilon_2$ of the equation $t^3+1=0$ have the property $\epsilon_1 \epsilon_2=1$,
it follows that $(\epsilon_i,1,0)= (1,\epsilon_j,0)$ in $\mathbb{P}^2_{\mathbb{K}}$, for $i\not = j$, $i,j=1,2$.
The diagram
$$
\begin{array}{ccc} \overline{DF} & (x=3at^2, y = 3at, z=1+t^3)&(x,y,z)\not= O\\ \
\bar{\bar p}\Big\uparrow \wr& \Big\uparrow&\Big\downarrow\\ \
A^1_\mathbb{K}&t&t=\frac{x}{y}\end{array}
$$
$$
\begin{array}{ccc}(0,0,1)&(1,-1,0)&(\epsilon_i,1,0)\\ \
\Big\downarrow \Big\uparrow&\Big\uparrow&\Big\uparrow\\ \
t=0&t=-1&t=\epsilon_i\end{array}
$$
shows that the parametrization $\bar {\bar p}$ is bijective.
The function $\bar {\bar p}$ transports the additive group structure from $\mathbb{K}$
to $\overline{DF}$ and we obtain the group $(\overline{DF},\oplus)$ defined by
$$\bar {\bar p}(t) \oplus  \bar{\bar p}(t^{\prime}) \,\,\overset{def}{=}\,\, \bar{\bar p}(t+t^{\prime})$$
or, equivalently, by
$$(3at^2,3at,1+t^3)\oplus (3a{t^{\prime}}^2, 3at^{\prime},1+{t^{\prime}}^3)\,\,\overset{def}{=}\,\, (3a(t+t^{\prime})^2,3a(t+{t^{\prime}}),1+(t+{t^{\prime}})^3).$$
It appears an isomorphism of groups
$$(\overline{DF}, \oplus)\stackrel{\bar {\bar p}}{\stackrel{\sim}{\longleftarrow}}(\mathbb{K}, +)$$
$$\bar{\bar p}(t)\longleftarrow t$$

The point $\bar{\bar p}(0) = (0,0,1)=O$ is a neutral element for the group $(\overline{DF}, \oplus)$.
Obviously,
$$\forall t \in \mathbb{K},\,\, -\bar{\bar p}(t)= \bar{\bar p}(-t).$$
In this way, $-\bar{\bar p}(-1)= \bar{\bar p}(1)$ implies the identification between
$-(1,-1,0)$ (the opposite of the point $I$ at infinite) and $V=(3a,3a,2)\in \overline{DF}$
(the vertex of $\overline{DF}\subset \mathbb{P}^2_{\mathbb{K}}$).

We have a natural isomorphism of groups:
$$\bar \sigma: (\overline{DF},+)\longrightarrow (\overline{DF},\oplus)$$
$$(x,y,z)\longrightarrow (y,x,z).$$

It is easy to check that: (1) the group composition laws $+$ and $\oplus$ on $\overline{DF}$ are distinct,
(2) for each $P\in \overline{DF}$, the symmetric/opposite $-P$ of $P$ is the same w.r.t.
each of the composition laws $+$ and $\oplus$.

If $\mathbb{K}=\mathbb{R}, \mathbb{C}$, the group $(\overline{DF},\oplus)$ is not a topological group if we consider $\overline{DF}$
with the topology induced by the real, resp. complex, topology of $\mathbb{P}^2_{\mathbb{R}}$,
resp. $\mathbb{P}^2_{\mathbb{C}}$.

Also, it is easy to see that the diagram
$$\begin{array}{ccc}(\overline{DF},+)& \underset{\sim}{\overset{\bar\sigma}\longrightarrow}&(\overline{DF},\oplus)\\ \
{\wr}{\nwarrow} \bar{p} & &{\wr}{\nearrow} \bar{\bar p} \\ \ &(\mathbb{K},+)& \end{array}$$
is commutative.

\section{Some topologies on projective Descartes Folium}

\hspace{0.5cm} Let $\mathbb{K}=\mathbb{R}, \mathbb{C}$.
Now, as for the previous first and second "exotic" structures, by using identical ideas,
we will show that the groups $(\overline{DF},+)$ and $(\overline{DF},\oplus)$ become
topological groups if we consider some finer topologies on $\overline{DF}$.

We have a similar situation as in the cases of the first and second "exotic" structures (see $\S$ 4.1 and $\S$ 4.2).

{\bf Proposition} {\it Let $\mathbb{K}=\mathbb{R}$ or $\mathbb{K}=\mathbb{C}$ and
$\overline{DF}\subset \mathbb{P}^2_{\mathbb{K}}$ endowed with the topology $\bar \tau$ induced by
the real, resp. complex, topology of $\mathbb{P}^2_{\mathbb{K}}$. Then:

(i) the bijective map $\bar p:\mathbb{A}^1_{\mathbb{K}}=\mathbb{K} \overset{\sim}\longrightarrow \overline{DF}$ is
continuous but not a homeomorphism;

(ii) the restriction
$\bar p:\mathbb{A}^1_{\mathbb{K}}\setminus \{O\}=\mathbb{K}\setminus\{0\} \overset{\sim}\longrightarrow \overline{DF}\setminus\{O\}$
is a homeomorphism.}

Similar properties hold for
$\bar{\bar p}=\bar\sigma\circ \bar p: \mathbb{A}^1_{\mathbb{K}}=\mathbb{K} \overset{\sim}\longrightarrow \overline{DF}$.

Recall that for each $t\in \mathbb{K}\setminus\{0\}$, we have
$\bar{\bar p}(t) ={\bar\sigma} \bar p(t) = {\bar p} \left(\frac{1}{t}\right)$.

We will introduce two topologies ${\bar\tau}_S$ and ${\bar\tau}_W$ on $\overline{DF}$ as follows:

{\bf Definitions} {\it ${\bar\tau}_S$ (resp. ${\bar\tau}_W$) is the image in $\overline{DF}$ of
the real, resp. complex, topology of $\mathbb{A}^1_{\mathbb{K}}=\mathbb{K}$ by the bijective map
$\bar p:\mathbb{K}\overset{\sim}\longrightarrow \overline{DF}$, resp.
$\bar{\bar p}={\bar\sigma}{\bar p}:\mathbb{K}\overset{\sim}\longrightarrow \overline{DF}$.}

It follows that $\bar \tau_S$ (resp. $\bar \tau_W$) is separated, connected, paracompact, locally compact and with countable basis.

Denote now
$$S =\bar p(\{z\in \mathbb{K}\,\Big\vert\, |z|<1\})\subset DF\subset \overline{DF}$$
$$W =\bar p(\{z\in \mathbb{K}\,\Big\vert\, |z|>1\})\cup\{O\}\subset DF\subset \overline{DF}.$$
We have
$$S\cap W=\{O\},\,\, S\cup W=\overline{DF}\setminus\{\bar p(z)\Big\vert z\in \mathbb{K}, |z|=1\}$$
$$\bar\sigma(S)=W,\,\, \bar\sigma(W)=S,\,\, W= \bar{\bar p}(\{z\in \mathbb{K}\Big\vert\,|z|<1\}).$$

For $\mathbb{K}=\mathbb{R}$, the branches $S$ and $W$ are just the subsets of
$DF$ defined previously for the case $DF \subset \mathbb{A}^2_{\mathbb{R}}$ (see $\S 4$.

The list of properties of $\tau_S$, resp. $\tau_W$ (see $\S$ 4.2), is valid also for
${\bar\tau}_S$, resp. ${\bar\tau}_W$, via the corresponding modifications. Let us recall them:

($i^\prime$) $\bar\tau_S$ (resp. $\bar\tau_W$) is a finer topology than $\bar \tau$ (i.e. $\bar \tau_S, \bar \tau_W\succ \bar \tau$).

($ii^\prime$) The induced topology ${\bar\tau}_S|_{\overline{DF}\setminus\{O\}}$
(resp. ${\bar\tau}_W|_{\overline{DF}\setminus\{O\}}$) on $\overline{DF}\setminus\{O\}\subset\mathbb{P}^2_{\mathbb{K}}$
is induced on $\overline{DF}\setminus\{O\}$ by the real, resp. complex, topology of $\mathbb{P}^2_{\mathbb{K}}$.

($iii^\prime$) If $\{U_n\}_{n\in \mathbb{N}}$ is a fundamental system of open neighborhoods of $O$ in $\mathbb{P}^2_{\mathbb{K}}$
w.r.t. real/complex topology, then $\{U_n\cap S\}_{n\in \mathbb{N}}$ (resp. $\{U_n\cap W\}_{n\in \mathbb{N}}$) is
a fundamental system of open neighborhoods of $O\in \overline{DF}$ w.r.t. the topology
$\bar \tau_S$ (resp. $\bar \tau_W)$.

($iv^\prime$) A basis for the topology $\bar \tau_S$ (resp. $\bar \tau_W$) is
$$\bar \tau \cup\{U\cap S\,|\,U \subseteq \mathbb{P}^2_{\mathbb{K}}\,\,\hbox{open subset}\}$$
$$\left(\hbox{resp.}\,\,\bar \tau \cup\{U\cap W\,|\,U \subseteq \mathbb{P}^2_{\mathbb{K}}\,\,\hbox{open subset}\}\right).$$
Moreover, for each $V \in \bar\tau_S$ (resp. $V\in \bar \tau_W$),
$$V = (U^\prime\cap \overline{DF})\cup(U\cap S)$$
$$\left( resp.\,\, V = (U^\prime\cap \overline{DF})\cup(U\cap W)\right),$$
with $U,U^\prime\subseteq \mathbb{P}^2_{\mathbb{K}}$ open subsets.

($v^\prime$) Let
$$\bar \pi:\overline{DF}\to \mathbb{K},\,\, \bar\pi= \left\{\begin{array}{ccc}\frac{y}{x}&if& (x,y,z)\not=O\\ \
0&if&(x,y,z)=O\end{array}\right.$$
$$\left(resp.\,\, \bar{\bar \pi}:\overline{DF}\to \mathbb{K},\,\, \bar{\bar\pi}= \left\{\begin{array}{ccc}\frac{x}{y}&if& (x,y,z)\not=O\\ \
0&if&(x,y,z)=O\end{array}\right.\right)$$
Then ${\bar \tau}_S$ (resp. ${\bar\tau}_W$) is the weakest topology on $\overline{DF}$ such that
$\bar \pi$ (resp.\,$\bar{\bar \pi}$) is continuous ($\mathbb{K}$ endowed with the real/complex topology).

($vi^\prime$) $\{O\}\subset W$ (resp.\, $\{O\}\subset S$) is a connected component
of the subspace $W$ (resp.\, $S$) w.r.t the topology ${\bar \tau}_S$ (resp. ${\bar\tau}_W$).
Moreover , if $\mathbb{K}=\mathbb{C}$, then
$$ W=\bar p(\{z\in \mathbb{K}\,\,\Big\vert\,\, |z|>1\})\cup\{O\}$$
$$ (\hbox{resp.}\,\,S=\bar{\bar p}(\{z\in \mathbb{K}\,\,\Big\vert\,\, |z|>1\})\cup\{O\})$$
is the representation of $W$ (resp. $S$) as the union of its connected components w.r.t. ${\bar\tau}_S$ (resp. ${\bar\tau}_W$).
On the other hand, $S$ (resp. $W$) is connected w.r.t. $\bar \tau_S$ (resp. $\bar \tau_W$).

\section{Some differential, resp. complex analytic structures on projective Descartes Folium}

\hspace{0.5cm} Let $\mathbb{K}=\mathbb{R},\mathbb{C}$. Let
$$\bar \pi: \overline{DF}\to \mathbb{K}\,\, (\hbox{resp.}\,\, \bar{\bar \pi}: \overline{DF}\to \mathbb{K})$$
defined in ($v^\prime$) above, where $\overline{DF}$ is endowed with the previous topology
${\bar\tau}_S$ (resp. ${\bar\tau}_W$). Then $\bar\pi = {\bar p}\,^{-1}$ (resp. $\bar{\bar\pi} = \bar{{\bar p}}\,^{-1}$)
is a homeomorphism.

Then the simple atlas $\{(\overline{DF},\bar \pi)\}$ (resp. $\{(\overline{DF},\bar{\bar \pi})\}$)
defines a structure of differential manifold (if $\mathbb{K}=\mathbb{R}$) or of complex analytic manifold
(if $\mathbb{K}=\mathbb{C}$) on the topological space $(\overline{DF},{\bar \tau}_S)$ (resp. $(\overline{DF},{\bar \tau}_W)$),
denoted by $\overline{\cal A}_S$ (resp. $\overline{\cal A}_Wq$).
Consequently $\bar \pi$ (resp. $\bar{\bar \pi}$) becomes a diffeomorphism (if $\mathbb{K}=\mathbb{R}$) or
an analytic isomorphism (if $\mathbb{K}=\mathbb{C}$).

Since
$$\bar p:(\mathbb{K},+)\overset{\sim}\longrightarrow (\overline{DF},+)$$
$$(\hbox{resp.}\,\, \bar{\bar p}:(\mathbb{K},+)\overset{\sim}\longrightarrow (\overline{DF},\oplus)$$
is also a group isomorphism, it follows

{\bf Proposition} {\it Let $\mathbb{K}=\mathbb{R},\mathbb{C}$. Then:

(i) $(\overline{DF},+)$ (resp. $(\overline{DF},\oplus)$) is a $\mathbb{K}$-Lie groups (in particular,
a topological group), where $\overline{DF}$ is endowed with the topology ${\bar\tau}_S$ (resp. ${\bar\tau}_W$)
and the differential or analytic manifold structure $\overline{\cal A}_S$ (resp. $\overline{\cal A}_W$)
given by the atlas $\{(\overline{DF},\bar\pi)\}$
(resp. $\{(\overline{DF},\bar{\bar\pi})\}$).

(ii)
$$\bar p:(\mathbb{K},+)\overset{\sim}\longrightarrow (\overline{DF},+)$$
$$(\hbox{resp.}\,\,\bar{\bar p}:(\mathbb{K},+)\longrightarrow (\overline{DF},\bar{\bar\pi}))$$
is an isomorphism of $\mathbb{K}$-Lie groups (in particular, of topological groups).}

It is obvious that in the commutative diagram
$$\begin{array}{ccccc}(\overline{DF},+)& &\underset{\sim}{\overset{\bar\sigma}\longrightarrow}&&(\overline{DF},\oplus)\\ \
&{\wr}{\nwarrow} \bar{p} & &{\wr}{\nearrow} \bar{\bar p}& \\ \ &&(\mathbb{K},+)&& \end{array}$$
$\bar \sigma =\bar{\bar p}\circ {\bar p}^{-1}$ is also an isomorphism of $\mathbb{K}$-Lie groups
(hence of topological groups).

\section{Projective Descartes Folium as \\topological field}

\hspace{0.5cm} Let $\mathbb{K}$ be a field with $\hbox{char.}\,\mathbb{K}\not = 3$ and $\overline{DF}\subset \mathbb{P}^2_{\mathbb{K}}$.

We already considered the commutative group $(\overline{DF}\setminus\{O\}, \cdot)$. The composition law
$\cdot$ can be extended trivially on whole $\overline{DF}$, defining
$$O\cdot A \,\,\overset{def}=\,\,O,\,\, A\cdot O \,\,\overset{def}=\,\,O$$
for each $A\in \overline{DF}$. Then
$$\bar p: (\overline{DF},\cdot) \overset{\sim}\longrightarrow (\mathbb{K}, \cdot)$$
and
$$\bar{\bar p}: (\overline{DF},\cdot) \overset{\sim}\longrightarrow (\mathbb{K}, \cdot)$$
are both isomorphisms of monoids. Then
$$\bar p: (\overline{DF},+,\cdot) \overset{\sim}\longrightarrow (\mathbb{K},+, \cdot)$$
and
$$\bar{\bar p}: (\overline{DF},\oplus,\cdot) \overset{\sim}\longrightarrow (\mathbb{K},+, \cdot)$$
are both bijective maps which are compatible with additive, resp. multiplicative composition laws.

Since $(\mathbb{K},+, \cdot)$ is a field, it follows that $(\overline{DF},+,\cdot)$ and $(\overline{DF},\oplus,\cdot)$
are both commutative fields. We have

{\bf Proposition} {\it Let $\mathbb{K}=\mathbb{R},\mathbb{C}$ and $\overline{DF}\subset \mathbb{P}^2_{\mathbb{K}}$.
Then:

(i) The field $(\overline{DF},+,\cdot)$  endowed with the topology ${\bar \tau}_S$,
and the field $(\overline{DF},\oplus,\cdot)$ endowed with the topology ${\bar \tau}_W$,
are both topological fields.

(ii) We have a commutative diagram of isomorphisms of topological fields
$$\begin{array}{ccccc}(\overline{DF},+,\cdot)&& \underset{\sim}{\overset{\bar\sigma}\longrightarrow}&&(\overline{DF},\oplus,\cdot)\\ \
&{\wr}{\nwarrow} \bar{p} & &{\wr}{\nearrow} \bar{\bar p}& \\ \ &&(\mathbb{K},+,\cdot)&& \end{array}$$}

In fact, $\bar p$ and $\bar{\bar p}$ are isomorphisms of fields and homeomorphisms. Since
$(\mathbb{K},+,\cdot)$ is a topological field, it follows that $(\overline{DF},+,\cdot)$
and $(\overline{DF},\oplus,\cdot)$ are both topological fields w.r.t. the topology ${\bar \tau}_S$,
resp. ${\bar \tau}_W$, and $\bar p,\,\,\bar{\bar p}$ are isomorphisms of topological fields.
The commutativity of the previous diagram is already known and consequently
$\bar \sigma= \bar{\bar p}\circ {\bar p}\,^{-1}$ is also an isomorphism of topological fields.

{\bf Acknowledgments}

Partially supported by "Simion Stoilow" Institute of Mathematics of the Romanian Academy,
University Politehnica of Bucharest, UNESCO Chair in Geodynamics-"Sabba S. \c Stef\u anescu"
Institute of Geodynamics of the Romanian Academy and by Academy of Romanian Scientists.

\framebox[1.1\width]{Adrian Constantinescu}, Simion Stoilow Institute of Mathematics of the Romanian Academy,
C.P. 1-764, RO-014700 Bucharest, Romania\\
Email: Adrian.Constantinescu@imar.ro\\

Prof. Dr. Constantin Udri¸ste, PhD student Steluta Pricopie, University Politehnica of Bucharest, Faculty of Applied
Sciences, Department of Mathematics-Informatics, Splaiul Independentei 313,
060042 Bucharest, Romania\\
E-mail: udriste@mathem.pub.ro, anet.udri@yahoo.com\\
E-mail: maty\_star@yahoo.com

\end{document}